\documentclass[11pt]{amsart}
\pdfoutput=1 

\usepackage[
text={440pt,575pt},
headheight=9pt,
centering
]{geometry}
\usepackage{appendix}
\usepackage{hyperref}
\usepackage{xcolor}

\definecolor{darkred}{RGB}{160,0,0}
\definecolor{darkblue}{RGB}{0,0,160}
\hypersetup{
  colorlinks,
  citecolor=darkblue,
  filecolor=black,
  linkcolor=darkblue,
  urlcolor=darkblue
}

\usepackage[utf8]{inputenc}
\usepackage[T1]{fontenc}
\usepackage{microtype}
\usepackage{lmodern}
\usepackage{amsmath,amssymb,amsthm}
\usepackage{soul}
\usepackage{mdwlist}
\usepackage[pdftex]{graphicx}
\usepackage{latexsym}
\usepackage{verbatim}
\usepackage{graphicx,caption,subcaption}
\usepackage{xcolor}
\usepackage{todonotes} 
\usepackage{enumitem}
\usepackage{upgreek}

\allowdisplaybreaks 

\theoremstyle{plain}
\newtheorem{theorem}{Theorem}[section]
\newtheorem{corollary}[theorem]{Corollary} 
\newtheorem{lemma}[theorem]{Lemma}

\theoremstyle{definition}
\newtheorem{definition}[theorem]{Definition} 	
\newtheorem{remark}[theorem]{Remark}

\newcommand{\R}{\mathbb{R}}

\newcommand{\N}{\mathbb{N}}
\newcommand{\PP}{\mathbb{P}}
\newcommand{\NN}{\mathfrak{N}}

\newcommand{\A}{\mathcal{A}}
\newcommand{\D}{\mathcal{D}}

\newcommand{\E}{\mathbb{E}}
\newcommand{\Var}{\mathrm{Var}}
\newcommand{\Cov}{\mathrm{Cov}}
\newcommand{\corr}{\mathrm{Corr}}

\newcommand{\inv}{\mathrm{inv}}
\newcommand{\des}{\mathrm{des}}

\title{Joint asymptotic normality and extremes of generalized inversions and descents 
}
\author[P. Dörr]{Philip Dörr}
\address{Department of Mathematics, Otto-von-Guericke University Magdeburg}
\email{philip.doerr@ovgu.de}

\date{}



\thanks{Philip Dörr is supported by the DFG (314838170, GRK 2297,
  ``MathCoRe'').}

\begin{document}

\begin{abstract}
Generalized inversions $X_\inv^{(d)}$ and generalized descents $X_\des^{(d)}$ are an interesting combinatorial extension of the common inversion and descent statistics on permutation groups. By means of the root poset, they can be defined on all classical Weyl groups (i.e., symmetric groups, signed and even-signed permutation groups). In this paper, we investigate the bivariate normality of $(X_\inv^{(d_1)}, X_\des^{(d_2)})^\top$ for different regimes $d_1, d_2$ that vary with the group size. Depending on the regimes, we classify the existence of the asymptotic correlation and provide explicit formulas. We also draw conclusions on the extreme value behavior of $X_\inv^{(d_1)}$, $X_\des^{(d_2)},$ and $(X_\inv^{(d_1)}, X_\des^{(d_2)})^\top$ in suitable triangular arrays, discussing bounds on the number of i.i.d.\ samples $k_n$ from a group of size $n$ so that the extreme values of these statistics are attracted to the Gumbel distribution.
\end{abstract}

\maketitle

\section{Introduction}

The numbers of inversions and descents are two important characteristics of permutations. For a permutation $\pi \negmedspace: \{1, \ldots, n\} \rightarrow \{1, \ldots, n\},$ an \textit{inversion} is a pair $(i,j)$ with $i < j,$ but $\pi(i) > \pi(j)$. A \textit{descent} is an index $i$ with $\pi(i) > \pi(i+1)$, i.e., descents correspond to adjacent inversions. Equipping the symmetric group $S_n$ with the discrete uniform probability measure induced by the point masses $\PP(\{\pi\}) = 1/n!$ $\forall \pi \in S_n$, we denote the random numbers of inversions and descents by $X_\inv$ and $X_\des$. These random variables can also be represented with help of i.i.d.\ variables $Z_1, Z_2, \ldots, \sim U(0,1)$ by means of
\begin{align}
    X_\inv &= \sum_{1\leq i<j\leq n} \textbf{1}\{Z_i > Z_j\}\,, \label{1} \\
    X_\des &= \sum_{i=1}^{n-1} \textbf{1}\{Z_i > Z_{i+1}\}\,. \nonumber 
\end{align}
In \eqref{1}, we are taking the sum of the indicators $\textbf{1}\{Z_i > Z_j\}$ over all pairs $(i,j)$ with $i<j$. A class of generalized inversion statistics can be constructed by restricting the sum in \eqref{1} to pairs $(i,j)$ with $1 \leq j-i \leq d$, for some $d \in \{1, \ldots, n-1\}$. Writing $\NN_{n,d} := \{(i,j) \in \{1,\ldots,n\}^2 \mid 1 < j-i \leq d\}$, we say that
\begin{equation} 
X_\inv^{(d)} := \sum_{(i,j) \in \NN_{n,d}} \textbf{1}\{Z_i>Z_j\} \label{3}
\end{equation}
counts the so-called $d$-\textit{inversions}. The term \textit{generalized inversions} is an umbrella term for all $d$-inversions. By choosing $d=n-1$ or $d=1,$ it is seen that this class includes common inversions and descents. 

Generalized inversions were first introduced by de Mari \& Shayman \cite{de1988generalized} to describe the Betti numbers of Hessenberg subvarieties in regular complex-valued matrices. The first probabilistic investigations of the random variables $X_\inv^{(d)}$ are due to Bona \cite{bona2007generalized} and Pike~\cite{pike2011convergence}, who computed the mean and variance. Bona proved a central limit theorem (CLT) for $X_\inv^{(d)}$ with $d$ remaining fixed, while Pike also addressed the case of $d = d(n)$ varying with $n$, and proved bounds for the rate of convergence depending on both $n$ and $d(n)$.  

Symmetric groups belong to the class of finite irreducible Coxeter groups, which have been classified by \cite{coxeter1935complete}. The concept of inversions and descents can be transferred to these groups, see \cite[Section 1.4]{bjorner2006combinatorics}. Two other important subfamilies of the finite irreducible Coxeter groups are the \textit{signed permutation groups} $B_n$ and the \textit{even-signed permutation groups} $D_n$. The group $B_n$ consists of all maps $\pi \negthickspace: \{1, \ldots, n\} \longrightarrow \{1, \ldots, n\} \cup \{-1, \ldots, -n\}$ for which $|\pi|$ is a permutation. Such a map is called a \textit{signed permutation}. The group $D_n$ is the subgroup of $B_n$ containing all signed permutations with an even number of negative signs. We refer to the three families $S_n, B_n,$ and $D_n$ as the \textit{classical Weyl groups}.

Recently, Meier \& Stump \cite{meier2022central} introduced \textit{generalized descents} alongside generalized inversions, and further extended the concept of both generalized inversions and generalized descents to the groups $B_n$ and $D_n$. On the symmetric groups, $d$-\textit{descents} or \textit{generalized descents} are given by
\begin{equation}
    X_\des^{(d)} := \sum_{i=1}^{n-d} \textbf{1}\{Z_i > Z_{i+d}\}\,. \label{4}
\end{equation}
Meier \& Stump extensively computed the mean and variance of $X_\inv^{(d)}$ and $X_\des^{(d)}$ on the groups $B_n$ and $D_n,$ and also stated CLTs for both $X_\inv^{(d)}$ and $X_\des^{(d)}$ on all classical Weyl groups. As $d = d(n)$ varies with $n,$ the CLTs tacitly require that the number of summands in \eqref{3} or \eqref{4} tends to infinity.

The joint dependency structure of $X_\inv^{(d_1)}$ and $X_\des^{(d_2)}$ for two different sequences $d_1 = d_1(n),$ $d_2 = d_2(n)$ is interesting. For common inversions and descents, the asymptotic correlation of $X_\inv$ and $X_\des$ is zero and the CLT for $(X_\inv, X_\des)^\top$ is satisfied, with the limit being the bivariate standard normal distribution. This was proven by Fang \& R\"ollin \cite{fang2015rates} for symmetric groups and recently by D\"orr \& Heiny \cite{dorr2025joint} for the other classical Weyl groups.

In this paper, we investigate the bivariate asymptotic normality of generalized inversions and descents. It turns out that the asymptotic correlation of $X_\inv^{(d_1)}$ and $X_\des^{(d_2)}$ can be strictly positive if both $d_1$ and $d_2$ grow linearly, i.e., if $d_1/n$ and $d_2/n$ do not tend to zero. In that case, the asymptotic correlation explicitly depends on the limits of $d_1/n$ and $d_2/n$. If these limits do not exist, then the limit of $\corr(X_\inv^{(d_1)}, X_\des^{(d_2)})$ commonly does not exist as well, hence, there is no CLT. Other interesting cases arise if $d_1$ and $d_2$ both remain bounded. Then, the asymptotic correlation turns out to be nonzero if $d_1 \geq d_2$ and if $d_1$ remains constant. In the extreme case of $d_1 = d_2 = 1,$ the two statistics $X_\inv^{(1)}$ and $X_\des^{(1)}$ are identical. If $(d_1(n))_{n \in \N}$ has multiple accumulation points, then there is no CLT. We will investigate this behavior in thorough detail.

Moreover, the investigation of extreme value asymptotics for common inversions and descents was initiated in \cite{dorr2022extreme, dorr2025joint}, so we aim to derive the corresponding results for generalized inversions and descents. As these permutation statistics are discrete-valued, extreme value limit theorems are formulated for triangular arrays of i.i.d.\ copies of $X_\inv^{(d_1)}$ and $X_\des^{(d_2)}$. The $n$-th row of such a triangular array corresponds to a classical Weyl group of rank $n$, and the number $k_n$ of i.i.d.\ copies is chosen depending on the rank. The row-wise maximum $M_n$ of these i.i.d.\ copies is considered and one seeks to find deterministic constants $a_n, b_n$ and an \textit{extreme value distribution} $G$ such that $a_n^{-1}(M_n - b_n) \overset{\D}{\longrightarrow} G$ holds. 

For permutation statistics that satisfy the CLT, the typical extreme value distribution $G$ is the Gumbel distribution, and the weak convergence $a_n^{-1}(M_n - b_n) \overset{\D}{\longrightarrow} G$ is simply called \textit{Gumbel attraction}. Typically, the extreme value behavior in triangular arrays is degenerate if there are too many i.i.d.\ samples per row (i.e., if $k_n$ is too large). We therefore seek to find asymptotic regimes of $k_n$ as large as possible so that Gumbel attraction of $X_\inv^{(d_1)},$ $X_\des^{(d_2)}$ and the joint statistic $(X_\inv^{(d_1)}, X_\des^{(d_2)})^\top$ is guaranteed. We thoroughly investigate the impact of $d_1$ and $d_2$ on these bounds for $k_n$ as well. 

This paper is structured as follows. Section~\ref{section2} reviews basic properties on generalized inversions and descents. Section~\ref{section3} is devoted to the study of the asymptotic correlation and the bivariate CLT of generalized inversions and descents. The CLT is achieved through high-dimensional Gaussian approximation rates and Haj\'{e}k projections.  Section~\ref{section4} addresses their extreme value limit theorems. The appendix sections~\ref{appendixB},~\ref{appendixC},~\ref{appendixA} collect all proofs. We use the shorthand notation $x \wedge y := \min\{x,y\}$  and we use typical Landau notation for positive sequences $a_n, b_n$ as follows: 
\begin{itemize}
    \item $a_n = O(b_n)$ means that 
    $\limsup_{n\rightarrow\infty} a_n/b_n < \infty$. 
    \item $a_n = o(b_n)$ means that   $\lim_{n\rightarrow\infty} a_n/b_n = 0$. This is also written as $b_n = \omega(a_n)$ or $b_n \gg a_n$.
    \item $a_n = \Theta(b_n)$ means that $a_n$ and $b_n$ have the same order of magnitude, i.e., both $a_n = O(b_n)$ and $b_n = O(a_n)$ hold.
    \item $a_n = b_n + o_\PP(1)$ means that $a_n, b_n$ are sequences of random variables with $a_n - b_n \overset{\PP}{\longrightarrow} 0$.
\end{itemize}

\section{Generalized inversions and descents} \label{section2} 

\begin{definition} \label{def2.1}
Let $S_n$ be a symmetric group and let $\pi \in S_n$. Recall $\NN_{n,d} := \{(i,j) \in \{1,\ldots,n\}^2: 1 < j-i \leq d\}$. For $d \in \{1, \ldots, n-1\}$,  $d$-\textit{inversions} are all pairs $(i,j)$ in $\NN_{n,d}$ with $\pi(i) > \pi(j)$. Moreover, $d$-\textit{descents} are all numbers $i \in \{1, \ldots, n-d\}$ with $\pi(i) > \pi(i+d)$. We write $X_\inv^{(d)}$ for the random number of $d$-inversions and $X_\des^{(d)}$ for the random number of $d$-descents. Probabilistic representations by means of i.i.d.\ variables $Z_1, Z_2, \ldots \sim U(0,1)$ are given in \eqref{3}, \eqref{4}. \footnote{In the literature, there are different terminologies, e.g., in \cite{bona2007generalized, pike2011convergence},  $d$-inversions are called $d$-descents. However, we use the terms and notation provided in \cite{meier2022central}.} 
\end{definition}

\begin{remark} \label{bem2.2}
Obviously, each $k \in \{1, \ldots, n\}$ can appear in at most $2d$ $d$-inversions. This bound is redundant if $d > n/2$. In fact, it is an important case distinction whether $d \leq n/2$ or $d > n/2$. In the former case, we can split $\{1, \ldots, n\}$ into the regions 
\[K_1 := \{1, \ldots, d\}, \qquad K_2 := \{d+1, \ldots, n-d\}, \qquad K_3 := \{n-d+1, \ldots, n\}\,,\] 
where $K_2 = \emptyset$ if $d=n/2$. For each $k \in K_2,$ \textit{all} larger indices $j = k+1, \ldots, k+d$ and \textit{all} smaller indices $i = k-1, \ldots, k-d$ allow to form $d$-inversions $(k, j)$ and $(i,k)$, respectively. For any $k \notin K_2,$ there are less than $d$ indices available in one direction. If $k \in K_1$, then only $k-1$ smaller indices $i < k$ are available for $d$-inversions, which is \textit{less} than $d$ indices. If $k > n-d$ is large, then there are only $n-k < d$ larger indices $j > k$ available. 

On the contrary, if $d > n/2,$ then $n-d < d$ and the above partition into three regions is written as 
\[K_1 := \{1, \ldots, n-d\}, \qquad K_2 := \{n-d+1, \ldots, d\}, \qquad K_3 := \{d+1, \ldots, n\}.\]
Now, if $k \in K_2,$ then all $k-1$ smaller indices $i<k$ and $n-k$ larger indices $j>k$ are available to form $d$-inversions. 
\end{remark}

\noindent The mean and variance of $X_\inv^{(d)}$ have been extensively computed by Pike \cite{pike2011convergence}. It is easy to verify that the special cases $d=1$ and $d=n-1$ are consistent with \cite[Corollaries 3.2~and~4.2]{kahle2020counting}. The proof of \cite[Theorem 1]{pike2011convergence} is reviewed in \cite[Theorem A.1]{meier2022central}, where the variance of $d$-descents is provided as well.

\begin{theorem}[see \cite{pike2011convergence}, Theorem 1 and \cite{meier2022central}, Theorem A.1] \label{thm2.3} 
For all $d = 1, \ldots, n-1,$ it holds that
\begin{align*}
    \E\left(X_\inv^{(d)}\right) &= \frac{2nd - d^2 - d}{4}, & \E\left(X_\des^{(d)}\right) &= \frac{n-d}{2}. 
\end{align*}
Moreover, if $d \leq n/2,$ then
\begin{align*}
    \Var\left(X_\inv^{(d)}\right) &= \frac{6nd + 4d^3 + 3d^2 - d}{72}, & \Var\left(X_\des^{(d)}\right) &= \frac{n+d}{12}.
\end{align*}
If $d > n/2,$ then
\[\Var\left(X_\inv^{(d)}\right) = -\frac{1}{6}d^3 + \left(\frac{1}{3}n - \frac{7}{24}\right)d^2 - \left(\frac{1}{6}n^2 - \frac{5}{12}n + \frac{1}{8}\right)d + \frac{1}{36}n^3 - \frac{1}{12}n^2 + \frac{1}{18}n\,\]
and $\Var\left(X_\des^{(d)}\right) = (n-d)/4$.
\end{theorem}

Meier \& Stump \cite{meier2022central} introduced an extension of generalized inversions and descents from symmetric groups to the other classical Weyl groups $B_n$ and $D_n$. This extension is based on the \textit{root poset} of a classical Weyl group. We refer to \cite[Section 2]{meier2022central} for explanations and details. On the symmetric group $S_n$, the ordered pairs of indices $(i,j)$ correspond to the positive roots $[ij] := e_i - e_j$, where $e_i, e_j$ are unit vectors in $\R^n$, and the height of $[ij]$ within the root poset is $\text{ht}([ij]) = j-i$.

On the signed permutation group $B_n,$ we also have to consider the positive roots $[\widetilde{ij}] := e_i + e_j$ and $[i] := e_i$ for $1 \leq i < j \leq n$. The heights of these additional roots are $\text{ht}([\widetilde{ij}]) = i+j$ and $\text{ht}([i]) = i$. See \cite[Example 2.2]{meier2022central} for an illustration of the root poset of $B_n$. On the even-signed permutation group $D_n,$ the roots $[i]$ are disregarded, and $[\widetilde{ij}]$ has height $i+j-2$. 

\begin{definition} \label{def2.3.6}
For any classical Weyl group, $d$-inversions are determined by roots of height at most $d$, and $d$-descents are determined by roots of height exactly $d$, see \cite[Definition 2.4]{meier2022central}. For symmetric groups, this coincides with Definition~\ref{def2.1}. In addition to $\NN_{n,d}$, we introduce 
\[
\widetilde{\NN}_{n,d} := \{(i,j) \in \{1,\ldots,n\}^2 \mid i<j, j+i \leq d\}\,, 
\]
with $d \in \{1, \ldots, 2n-1\}.$ Then, on the signed and even-signed permutation groups, $X_\inv^{(d)}$ and $X_\des^{(d)}$ can be expressed as follows: 
\begin{subequations}
\begin{align}
    X_\inv^{B,(d)} &= \sum_{(i,j) \in \NN_{n,d}} \textbf{1}\{Z_i > Z_j\} + \sum_{(i,j) \in \widetilde{\NN}_{n,d}} \textbf{1}\{-Z_i > Z_j\} + \sum_{i=1}^{n \wedge d} \textbf{1}\{Z_i < 0\}\,, \label{2.7a} \\[1ex]
    X_\inv^{D,(d)} &= \sum_{(i,j) \in \NN_{n,d}} \textbf{1}\{Z_i > Z_j\} + \sum_{(i,j) \in \widetilde{\NN}_{n,d+2}} \textbf{1}\{-Z_i > Z_j\}\,, \label{2.7b} \\[1ex]
    X_\des^{B,(d)} &= \sum_{i=1}^{n-d} \textbf{1}\{Z_i > Z_{i+d}\} + \sum_{i=1}^{\lceil d/2\rceil - 1} \textbf{1}\{-Z_i > Z_{d-i}\} + \textbf{1}\{Z_d < 0\}\,, \label{2.7c} \\[1ex]
    X_\des^{D,(d)} &= \sum_{i=1}^{n-d} \textbf{1}\{Z_i > Z_{i+d}\} + \sum_{i=1}^{\lceil d/2\rceil} \textbf{1}\{-Z_i > Z_{d+2-i}\}\,. \label{2.7d}
\end{align}    
\end{subequations}
Note that in \eqref{2.7c} and \eqref{2.7d}, indicators are ignored if they involve indices out of bounds. The largest possible choice of $d$ equals the total height of the root poset, namely $d_{\max} - 1$, where $d_{\max}$ denotes the largest degree of the underlying classical Weyl group. In particular, $d_{\max} = n$ for $S_n,$ $d_{\max} = 2n$ for $B_n$, and $d_{\max} = 2n-2$ for $D_n$. To precisely compute the variance of $X_\inv^{(d)}$ and $X_\des^{(d)}$ on the groups $B_n$ and $D_n,$ one needs to distinguish eight cases, as seen in \cite[Theorems A.4 and A.13]{meier2022central}. However, many of these cases give the same asymptotic quantification, which can be stated as follows:  
\end{definition}

\begin{lemma}[cf. \cite{meier2022central}, Theorems A.4 and A.13] \label{lemma2.3.8}
For the generalized inversions and descents on both the groups $B_n$ and $D_n,$ it holds that
\begin{align*}
    \Var\left(X_\inv^{(d)}\right) &= \begin{cases} \frac{1}{36}d^3 + \frac{1}{12}nd + O(d^2), & d \leq n/2 \\ \frac{1}{36}d^3 + O(d^2), & n/2 \leq d < n \\ -\frac{1}{12}d^3 + \frac{1}{3}nd^2 - \frac{1}{3}n^2d + \frac{1}{9}n^3 + O(d^2), & d \geq n \end{cases}, \\
    \Var\left(X_\des^{(d)}\right) &= \begin{cases}
        \frac{1}{24}d + \frac{1}{12}n + O(1), & d < n \\
        -\frac{1}{8}d + \frac{1}{4}n + O(1), & d \geq n
    \end{cases}.
\end{align*}
\end{lemma}

\section{Bivariate Asymptotic Normality} \label{section3}

The main work of this paper is to prove and investigate the asymptotic normality of the joint permutation statistic $\left(X_\inv^{(d_1)}, X_\des^{(d_2)}\right)^\top$, where $d_1 = d_1(n),$ $d_2 = d_2(n)$ can depend on the group size $n$. All proofs in this section are deferred to the appendix.

We will find that for low regimes of $d_1,$ the dependency structure of $\left(X_\inv^{(d_1)}, X_\des^{(d_2)}\right)^\top$ is sufficiently sparse, so that $\left(X_\inv^{(d_1)}, X_\des^{(d_2)}\right)^\top$ can be tackled by Gaussian approximation. For larger regimes, the complexity of the dependency structure needs to be reduced, and we aim to do so by replacing $X_\inv^{(d_1)}$ with its H\'{a}jek projection. In Section~\ref{sec3.1}, we introduce high-dimensional Gaussian approximation and prove that $\left(\hat{X}_\inv^{(d_1)}, X_\des^{(d_2)}\right)^\top$ exhibits the same asymptotic behavior as  $\left(X_\inv^{(d_1)}, X_\des^{(d_2)}\right)^\top$. Subsequently, in Section~\ref{sec3.2} we thoroughly investigate the covariance of $\hat{X}_\inv^{(d_1)}$ and $X_\des^{(d_2)}$ and the limit behavior of the correlation $\corr\left(\hat{X}_\inv^{(d_1)}, X_\des^{(d_2)}\right) = \Cov\left(\hat{X}_\inv^{(d_1)}, X_\des^{(d_2)}\right)/\left(\sigma(\hat{X}_\inv^{(d_1)})\sigma(X_\des^{(d_2)})\right)$.
In Section~\ref{sec3.3}, we investigate the asymptotic behavior of the correlation of the original joint statistic  $\left(X_\inv^{(d_1)}, X_\des^{(d_2)}\right)^\top$ for the low regimes, where the use of the H\'{a}jek projection is not feasible. 

\subsection{Gaussian approximation for sparse dependency structures} \label{sec3.1} 
The classical multivariate CLT states that for a sequence or a triangular array $(X_n)_{n \in \N} \subseteq \R^p$ and an appropriate $p$-variate normal distribution $\mathcal{N}$, under suitable conditions, the Gaussian approximation error
\begin{equation}
    \sup_{\textbf{x} \in \R^p} |\PP(X_n \leq \textbf{x}) - \PP(\mathcal{N} \leq \textbf{x})| = \sup_{A = (-\infty,\textbf{x}]^p} |\PP(X_n \in A) - \PP(\mathcal{N} \in A)| \label{gausserror}
\end{equation}
tends to zero as $n \rightarrow \infty$. Even stronger statements are obtained if the system $\A = \{(-\infty,\textbf{x}]^p\}$ of negative orthants is extended to, e.g., the class of all hyperrectangles or all convex sets, and if the dimension $p$ of the samples $X_1, \ldots, X_n$ can grow with $n$. The latter is commonly referred to as the \textit{high-dimensional setting} and it encompasses the fixed-dimensional setting by artificially repeating components of a random vector in fixed dimension.

Seminal works for high-dimensional Gaussian approximation involving independent random vectors are due to \cite{chernozhukov2013gaussian, chernozhukov2017central}, followed by several efforts to reduce the rate of convergence in \eqref{gausserror} and the growth rate of $p$, see \cite{chernozhukov2023nearly, das2021central, fang2021high, koike2021notes, kuchibhotla2020high}. The uniform consideration of all hyperrectangles yields not only CLTs but also extreme value limit results, which will be used in Section~\ref{section4}. 

An even more significant extension of Gaussian approximation is the consideration of more general dependency structures within the random vectors $X_1, \ldots, X_n,$ dropping the assumption of independence. A seminal work in this direction was achieved by Chang \textit{et al.}\ \cite{chang2024central}, which classifies the dependence of random vectors by dependency graphs. Let $\left(X_j^{(n)}\right)_{j=1, \ldots, n}$ be a triangular array of $p$-dimensional centered random vectors (without restriction), and let 
\begin{align*}
X^{(n)} &:= \sum_{t=1}^n X_t^{(n)}, & \Sigma^{(n)} &:= \Var(X^{(n)})\,, & \mathcal{N}_n &\sim \mathrm{N}_p(0,\Sigma^{(n)})\,.
\end{align*}
Let $\A^{\mathrm{re}}$ be the system of all $p$-dimensional hyperrectangles including infinite bounds, i.e., $\A^{\mathrm{re}} := \Bigl\{\{\textbf{w} \in \R^p \negmedspace: \textbf{a} \leq \textbf{w} \leq \textbf{b}\} \mid \textbf{a}, \textbf{b} \in [-\infty, \infty]^p\Bigr\}\,.$ Let 
\[
r_n(\A^{\mathrm{re}}) := \sup_{A \in \A^{\mathrm{re}}} |\PP(X^{(n)} \in A) - \PP(\mathcal{N}_n \in A)|\,.
\]
In particular, $r_n(\A^{\mathrm{re}})$ is an upper bound for 
\[
    r_n(\A^{\mathrm{CLT}}) := \sup_{u \in \R^p} |\PP(X^{(n)} \leq u) - \PP(\mathcal{N}_n \leq u)|\,,
\]
which is relevant for the CLT, and for
\[
    r_n(\A^{\mathrm{ext}}) := \sup_{u \in \R^p} |\PP(X^{(n)} > u) - \PP(\mathcal{N}_n > u)|\,,
\]
which is relevant for the asymptotics of extreme values. Now, let $G_n$ be a dependency graph for $X_1^{(n)}, \ldots, X_n^{(n)},$ which consists of all edges $(i,j)$ for which $X_i^{(n)}$ and $X_j^{(n)}$ are dependent. Let $\Delta_n$ be the maximum degree of $G_n$ and let $\Delta_n^*$ be the maximum degree of the 2-reachability graph of $G_n$. If the graphs $G_n$ are not too dense as $n \rightarrow \infty,$ then $r_n(\A^{\mathrm{re}})$ can be bounded as follows:

\begin{theorem}[cf. \cite{chang2024central}, Theorem 2] \label{thm3.3}
Let $(X_j^{(n)})_{j=1,\ldots,n}$ be a triangular array and let $r_n(\A^{\mathrm{re}}),$ $\Delta_n,$ and $\Delta_n^*$ be as above. Under some regularity conditions which are satisfied by bounded and non-degenerate random variables, it holds that
\[
    r_n(\A^{\mathrm{re}}) = O\left(\frac{(\Delta_n\Delta_n^*)^{1/3}\log(p)^{7/6}}{n^{1/6}}\right).
\]
In particular, if $\Bigl(X_j^{(n)}\Bigr)_{j=1,\ldots,n}$ is $m$-dependent for some global constant $m \in \N,$ then $r_n(\A^{\mathrm{re}}) = O\bigl(n^{-1/6}\log(p)^{7/6}\bigr)$. 
\end{theorem}

In what follows, we consider the joint statistic $\left(X_\inv^{(d_1)}, X_\des^{(d_2)}\right)^\top$ for two sequences $d_1 = d_1(n),$ $d_2 = d_2(n)$, with $d_1(n), d_2(n) \leq d_{\max} - 1$ $\forall n \in \N$. Both $X_\inv^{(d_1)}$ and $X_\des^{(d_2)}$ are based on a sequence of classical Weyl groups $(W_n)_{n \in \N},$ where each $W_n$ is one of $S_n, B_n,$ or $D_n$. We keep using $d$ as an umbrella notation for $d_1$ or $d_2,$ depending on the context.

If $d_2(n)$ diverges as $n\rightarrow\infty$, then there is no $m \in \N$ for which all $X_\des^{(d)}$ are $m$-dependent, but the dependency structure of $X_\des^{(d)}$ is still sparse. Recall that $X_\inv^{(d)}$ and $X_\des^{(d)}$ are constructed from i.i.d.\ random variables $Z_1, \ldots, Z_n$. According to \eqref{4}, \eqref{2.7c}, \eqref{2.7d}, we can represent $X_\des^{(d)}$ as a sum of indicator variables, each of which depends on at most three others. Therefore, the first-order and second-order maximum degrees $\Delta_n, \Delta_n^*$ are bounded in the way of $\Delta_n \leq 3$ and $\Delta_n^* \leq 9$. 
However, it has to be ensured that $X_\des^{(d)}$ consists of a growing number of summands, which was already noted in \cite[Corollary 2.7]{meier2022central}. This number of summands equals the number of positive roots with height exactly $d$, which we denote by $N_{n,d}^=$. \par 
The maximum degrees $\Delta_n, \Delta_n^*$ of the dependency graphs of the joint statistic $\left(X_\inv^{(d_1)}, X_\des^{(d_2)}\right)^\top$ are bounded in the way of $\Delta_n \leq 4d_1$ and $\Delta_n^* \leq 8d_1$. Moreover, we have to take into account that by \eqref{3}, \eqref{2.7a}, \eqref{2.7b}, $X_\inv^{(d_1)}$ is based on $\Theta(nd_1)$ summands, so we have to replace $n$ with $nd$ in Theorem~\ref{thm3.3}. Therefore, Theorem~\ref{thm3.3} gives an $o(1)$ bound of $r_n(\A^{\mathrm{re}})$ only if $d_1^{2/3}\log(nd_1)^{7/6} = o((nd_1)^{1/6})$, which leads to the condition $d_1 = o\left(n^{1/3}\log(n^{4/3})^{-7/3}\right)$. The consideration of these \textit{low regimes} is postponed to Section~\ref{sec3.3}. 

For any larger regime of $d_1,$ the dependency structure of $\left(X_\inv^{(d_1)}, X_\des^{(d_2)}\right)^\top$ is too complex to apply Theorem~\ref{thm3.3}. For these regimes, we resort to H\'{a}jek projections in order to approximate $X_\inv^{(d)}$ with a sum of independent variables. 

\begin{definition} 
Let $Z_1, \ldots, Z_n$ be independent random variables, and let $X$ be another random variable. Then, the \textit{H\'{a}jek projection} of $X$ with respect to $Z_1, \ldots, Z_n$ is given by
\[\hat{X} := \sum_{k=1}^n \E(X \mid Z_k) - (n-1)\E(X).\]
\end{definition}
As every $\E(X \mid Z_k)$ is a measurable function only in $Z_k$, the H\'{a}jek projection is a sum of independent random variables. If $(X_n)_{n \in \N}$ is a sequence of random variables with each $X_n$ depending on $Z_1, \ldots, Z_n$, then the H\'{a}jek projections $\hat{X}_n$ should give a sufficiently accurate approximation to $X_n$. This can be ensured by the following criterion.

\begin{theorem}[cf.~\cite{van2000asymptotic}, Theorem 11.2] \label{thm2.8}
Consider a sequence $(X_n)_{n\ge 1}$ of random variables and their associated H\'{a}jek projections $(\hat{X}_n)_{n\ge 1}$. If $\Var(\hat{X}_n) \sim \Var(X_n)$ as $n \rightarrow \infty,$ then
\[\frac{X_n - \E(X_n)}{\Var(X_n)^{1/2}} = \frac{\hat{X}_n - \E(\hat{X}_n)}{\Var(\hat{X}_n)^{1/2}} + o_{\PP}(1).\]
\end{theorem} 

In the case of $d_1 = \Omega\left(n^{1/3}\log(n^{4/3})^{-7/3}\right)$, we have to replace $X_\inv^{(d_1)}$ with the H\'{a}jek projection $\hat{X}_\inv^{(d_1)}$, which is based on the representations \eqref{3}, \eqref{2.7a}, and \eqref{2.7b}. The validity of the condition in Theorem~\ref{thm2.8} can be classified once we compute the variance of $\hat{X}_\inv^{(d_1)}$.

\begin{lemma} \label{lemma3.4}
For the H\'{a}jek projection $\hat{X}_\inv^{(d)}$ of generalized inversions on symmetric groups, it holds that
\[
    \Var\left(\hat{X}_\inv^{(d)}\right) = \frac{1}{72}\left(4d^3 + 6d^2 + 2d\right)
\]
if $d \leq n/2,$ and
\[
    \Var\left(\hat{X}_\inv^{(d)}\right) = -\frac{1}{6}d^3 + \left(\frac{1}{3}n - \frac{1}{4}\right)d^2 - \left(\frac{1}{6}n^2 - \frac{1}{3}n + \frac{1}{12}\right)d + \frac{1}{36}n^3 - \frac{1}{12}n^2 + \frac{1}{18}n
\]
if $d > n/2$.
\end{lemma}

\begin{lemma} \label{lemma3.7}
For the H\'{a}jek projection $\hat{X}_\inv^{(d)}$ of generalized inversions on both the signed permutation groups $B_n$ and the even-signed permutation groups $D_n,$ it holds that
\[
    \Var\left(\hat{X}_\inv^{(d)}\right) = \begin{cases} d^3/36 + O(d^2)\,, & d \leq n \\ -d^3/12 + nd^2/3 - n^2d/3 + n^3/9\,, & d > n \end{cases}.
\]
\end{lemma}

By recalling the asymptotic quantifications of $\Var(X_\inv^{(d)})$ given in Theorem~\ref{thm2.3} and Lemma~\ref{lemma2.3.8}, we obtain the following conclusion:

\begin{corollary} \label{cor3.5}
Consider the generalized inversion statistic $X_\inv^{(d)}$ on classical Weyl groups $(W_n)_{n \in \N},$ where each $W_n$ is one of $S_n, B_n,$ or $D_n,$ and $d = d(n)$ satisfies $1 \leq d \leq d_{\max}-1$ $\forall n \in \N$. 
Then, $\Var\left(\hat{X}_\inv^{(d)}\right) \sim \Var\left(X_\inv^{(d)}\right)$ holds if and only if $d \gg \sqrt{n}$.
\end{corollary}

From Theorem~\ref{thm3.3}, Corollary~\ref{cor3.5} and Slutsky's lemma, we can conclude that $(X_\inv^{(d_1)}, X_\inv^{(d_2)})$ satisfies the CLT, provided that $N_{n,d_2}^= \longrightarrow \infty$, and $d_1$ is in the aforementioned regimes $d_1 = o\left(n^{1/3}\log(n^{4/3})^{-7/3}\right)$ or $d_1 = \omega(n^{1/2})$, and the correlation $\hat{\rho}_n$ of $\hat{X}_\inv^{(d_1)}$ and $X_\des^{(d_2)}$ has a limit, which is the same as the limit of $\rho_n = \corr(X_\inv^{(d_1)}, X_\des^{(d_2)})$. 

\begin{corollary}
If either $d_1 = o\left(n^{1/3}\log(n^{4/3})^{-7/3}\right)$ or $d_1 \gg \sqrt{n},$ and if $d_2 = d_2(n) \in \{1, \ldots, d_{\max}-1\}$ is chosen so that $N_{n,d_2}^= \longrightarrow \infty,$ and if
\[
    \rho := \lim_{n\rightarrow \infty} \hat{\rho}_n := \lim_{n \rightarrow \infty} \corr(\hat{X}_\inv^{(d_1)}, X_\des^{(d_2)})
\]
exists, then $\left(X_\inv^{(d_1)}, X_\des^{(d_2)}\right)^\top$ satisfies the CLT, that is,
\[
    \left(\frac{X_\inv^{(d_1)} - \E(X_\inv^{(d_1)})}{\sqrt{\Var(X_\inv^{(d_1)})}},\, \frac{X_\des^{(d_2)} - \E(X_\des^{(d_2)})}{\sqrt{\Var(X_\des^{(d_2)})}}\right)^\top \overset{\D}{\longrightarrow} \mathcal{N}_2\left(0, I_\rho\right), \qquad I_\rho := \begin{pmatrix} 1 & \rho \\ \rho & 1 \end{pmatrix}.
\]
\end{corollary}

\subsection{Asymptotic correlation for large regimes of \texorpdfstring{$d_1$}{d1}} \label{sec3.2}

We now investigate and classify the asymptotic correlation $\rho$ of $\hat{X}_\inv^{(d_1)}$ and $X_\des^{(d_2)}$ for all regimes $d_1 \gg \sqrt{n}$. 

\subsection*{Symmetric groups} 
Roughly speaking, the asymptotic correlation $\rho$ is almost always zero if $d_1$ or $d_2$ are sublinear, or if the remainder $n - d_2$ is sublinear. If both $d_1$ and $d_2$ grow linearly in the way of $d_1/n \rightarrow c_1, d_2/n \rightarrow c_2$ with $c_1, c_2 \in (0,1),$ then $\rho$ exists but explicitly depends on $c_1$ and $c_2$. If the sequences $d_1/n, d_2/n$ have multiple accumulation points, then the CLT almost never holds. 

\begin{lemma} \label{lemma3.10}
For symmetric groups, the covariance of the H\'{a}jek projection $\hat{X}_\inv^{(d_1)}$ and $X_\des^{(d_2)}$ is given as follows: If $d_1 \leq d_2,$ then
\begin{align*}
    \Cov\left(\hat{X}_\inv^{(d_1)}, X_\des^{(d_2)}\right) &= \begin{cases} d_1(d_1 + 1)/12, & d_1 + d_2 \leq n \\ n(d_1 + d_2)/6 - d_1d_2/6 - d_2(d_2+1)/12 - n(n-1)/12, & d_1 + d_2 > n \end{cases}.
\end{align*}
In particular, if $d_1 < d_2$, then $\Cov\left(\hat{X}_\inv^{(d_1)}, X_\des^{(d_2)}\right) = \Cov\left(X_\inv^{(d_1)}, X_\des^{(d_2)}\right)$. If $d_1 > d_2,$ then
\[
    \Cov\left(\hat{X}_\inv^{(d_1)}, X_\des^{(d_2)}\right) = \begin{cases} d_1d_2/6 - d_2(d_2-1)/12, & d_1 + d_2 \leq n \\ n(d_1+d_2)/6 - d_2^2/6 - d_1(d_1+1)/12 - n(n-1)/12, & d_1 + d_2 > n \end{cases}.
\]
\end{lemma} 

We can now precisely characterize the situations in which the limiting correlation $\rho$ exists, and compute $\rho$ accordingly.

\begin{theorem} \label{thm3.10}
On the symmetric groups, the limiting correlation 
\[
    \rho = \lim_{n \rightarrow \infty} \corr\left(\hat{X}_\inv^{(d_1)}, X_\des^{(d_2)}\right)
\]
is classified as follows. First, assume $d_1 = o(n),$ while $d_1 = \omega(n^{1/2})$.
\begin{enumerate}
    \item[\emph{a)}] If either $n - d_2 \ll d_1$ or $n - d_2 \gg d_1,$ then $\rho = 0$.
    \item[\emph{b)}] If $d_2 \asymp n - cd_1$ for some $c \geq 1,$ then $\rho = 1/(\sqrt{2}c)$. 
    \item[\emph{c)}] If $d_2 \asymp n - cd_1$ for some $c \in (0, 1),$ then $\rho = \sqrt{2c} (1 - c/2)$. 
\end{enumerate}
In what follows, assume $d_1/n \rightarrow c_1$ and $d_2/n \rightarrow c_2$ with $0 < c_1 \leq 1$ and $0 \leq c_2 \leq 1$. Put 
\begin{align*}  
    p(c_1) &:= -c_1^3/2 + c_1^2 - c_1/2 + 1/12\,, \\ 
    f(c_1, c_2) &:= \sqrt{1 - c_2}(2c_1 + c_2 - 1)\,, \\
    g(c_1, c_2) &:= (2c_1c_2 - c_2^2)/\sqrt{1 + c_2}\,, \\
    h(c_1, c_2) &:= (2c_2(1 - c_2) - (1 - c_1)^2)/\sqrt{p(c_1)}\,.
\end{align*}
Then:
\begin{enumerate}
    \item[\emph{d)}] If $c_1 \leq c_2 \leq 1/2,$ then $\rho = \sqrt{6}/2 \cdot \sqrt{c_1/(1+c_2)}.$
    \item[\emph{e)}] If $c_1 \leq 1/2$ and $c_2 > 1/2,$ then
    \[
        \rho = \begin{cases} \sqrt{c_1}/\sqrt{2(1 - c_2)}\,, & c_2 \leq 1 - c_1 \\ f(c_1, c_2)/(\sqrt{2}c_1^{3/2})\,, & c_2 > 1 - c_1  \end{cases}.
    \]
    \item[\emph{f)}] If $1/2 < c_1 \leq c_2,$ then $\rho = \sqrt{3}/6 \cdot f(c_1, c_2)/\sqrt{p(c_1)}$.
    \item[\emph{g)}] If $1/2 \geq c_1 > c_2 > 0,$ then $\rho = \sqrt{6}g(c_1, c_2)/(2c_1^{3/2})$.
    \item[\emph{h)}] If $c_1 > 1/2$ and $c_2 \leq 1 - c_1,$ then $  \rho = g(c_1, c_2)/(2\sqrt{p(c_1)}).$
    \item[\emph{i)}] If $c_1 > c_2 > 1 - c_1,$ then
    \[
        \rho = \begin{cases} h(c_1, c_2)/(2\sqrt{1 + c_2})\,, & c_2 \leq 1/2 \\ \sqrt{3}/6 \cdot h(c_1, c_2)/\sqrt{1 - c_2}\,, & c_2 > 1/2 \end{cases}.
    \]
\end{enumerate} 
In particular, if $c_2 = 0,$ then always $\rho = 0$. If $c_2 = 1$ and $c_1 > 0,$ then we have $\rho = 0$ as well. 
\end{theorem}

\subsection*{Signed and even-signed permutation groups}

The groups $B_n$ and $D_n$ do not differ in terms of the asymptotic behavior of $\Cov\left(\hat{X}_\inv^{(d_1)}, X_\des^{(d_2)}\right)$ or the asymptotic correlation of $\hat{X}_\inv^{(d_1)}$ and $X_\des^{(d_2)}$. Therefore, we list the magnitudes of $\Cov\left(\hat{X}_\inv^{(d_1)}, X_\des^{(d_2)}\right)$ and the asymptotic correlations only for the signed permutation groups $B_n$.

\begin{lemma} \label{lemma3.12}
Consider the generalized inversion and descent statistics $X_\inv^{B,(d_1)},$ $X_\des^{B,(d_2)}$ on the signed permutation groups. Up to an error in the scale of $O(d_1),$ the covariance of the H\'{a}jek projection $\hat{X}_\inv^{B,(d_1)}$ and $X_\des^{B,(d_2)}$ is given as follows: If $d_1 \leq d_2,$ then
\[
    \Cov\left(\hat{X}_\inv^{B,(d_1)}, X_\des^{B,(d_2)}\right) = \begin{cases} d_1(d_1+1)/24, & d_1 + d_2 \leq 2n \\ (2n - d_2 + 1)(d_1 + d_2/2 - n)/12, & d_1 + d_2 > 2n \end{cases}.
\]
If $d_1 > d_2$ and $d_1 + d_2 \leq 2n,$ then
\[
    \Cov\left(\hat{X}_\inv^{B,(d_1)}, X_\des^{B,(d_2)}\right) = \frac{1}{24}d_2(d_2 + 1) + \frac{1}{12}d_2(d_1 - d_2)\,.
\]
If $d_1 > d_2$ and $d_1 + d_2 > 2n,$ then in case of $d_2 < n,$
\[
    \Cov\left(\hat{X}_\inv^{B,(d_1)}, X_\des^{B,(d_2)}\right) = -\frac{1}{24}d_1(d_1 - 4n + 3) - \frac{1}{12}d_2(d_2 - 2n + 1) - \frac{n^2}{6} + \frac{n}{4}\,,
\]
and in case of $d_2 > n,$
\[
    \Cov\left(\hat{X}_\inv^{B,(d_1)}, X_\des^{B,(d_2)}\right) = -\frac{1}{24}d_1(d_1 - 4n + 3) - \frac{1}{12}d_2(d_2 - 2n - 1) - \frac{n^2}{6} + \frac{n}{12}\,.
\]
\end{lemma}

We now make observations similar to those for symmetric groups regarding the limiting correlation $\rho$ appearing in the CLT for $(X_\inv^{(B,d_1)}, X_\des^{(B,d_2)})$ and $(X_\inv^{(D,d_1)}, X_\des^{(D,d_2)})$. The proof of the next theorem works analogously to the proof of Theorem~\ref{thm3.10} and is therefore omitted.

\begin{theorem} \label{thm3.12}
On the signed and even-signed permutation groups, the limiting correlation of $\hat{X}_\inv^{(d_1)}$ and $X_\des^{(d_2)}$ is classified as follows. First, assume $d_1 = o(n),$ while $d_1 = \omega(n^{1/2})$. In that case, we observe the same results as in Theorem \emph{\ref{thm3.10},} i.e.,
\begin{itemize}
    \item If either $n - d_2 \ll d_1$ or $n - d_2 \gg d_1,$ then $\rho = 0$.
    \item If $d_2 \asymp n - cd_1$ for some $c \geq 1,$ then $\rho = 1/(\sqrt{2}c)$. 
    \item If $d_2 \asymp n - cd_1$ for some $c \in (0, 1),$ then $\rho = \sqrt{2c}(1 - c/2)$. 
\end{itemize}
Now, assume $d_1/n \rightarrow c_1$ and $d_2/n \rightarrow c_2$ for some $c_1, c_2 \in (0, 2]$. Put
\begin{align*}
    p(c_1) &:= -\frac{1}{12}c_1^3 + \frac{1}{3}c_1^2 - \frac{1}{3}c_1 + \frac{1}{9}\,, &
    f(c_1, c_2) &:= \frac{1}{4}c_1(c_1 - 4) + \frac{1}{2}c_2(c_2 - 2) + 1\,.
\end{align*}
Then:
\begin{itemize}
\item If $c_1 \leq c_2 \leq 1,$ then 
\[
    \rho = \frac{\sqrt{3}}{2} \cdot \frac{\sqrt{c_1}}{\sqrt{1 + c_2/2}}\,.
\]
\item If $c_2 < c_1 \leq 1,$ then 
\[
    \rho = \sqrt{3}\,\frac{c_2(c_1 - c_2/2)}{c_1^{3/2}\sqrt{1+c_2/2}}\,. 
\]
\item If $c_1 \leq 1,$ $c_2 > 1,$ and $c_1 + c_2 \leq 2,$ then 
\[
    \rho = \frac{1}{\sqrt{2}} \frac{\sqrt{c_1}}{\sqrt{2-c_2}}\,.
\]
\item If $c_1 \leq 1, c_2 > 1,$ and $c_1 + c_2 > 2,$ then
\[
    \rho = \sqrt{2}\,\frac{\sqrt{2 - c_2}}{c_1^{3/2}}\left(c_1 + \frac{c_2}{2} - 1\right).
\]
\item If $1 < c_1 \leq c_2 \leq 2,$ then
\[
    \rho = \frac{1}{\sqrt{6}}\,\frac{\sqrt{2-c_2}}{\sqrt{p(c_1)}}\left(c_1 + \frac{c_2}{2} - 1\right).
\]
\item If $2 \geq c_1 > c_2 > 1,$ then
\[
    \rho = -\frac{\sqrt{6}}{3} \frac{f(c_1, c_2)}{\sqrt{(2 - c_2)p(c_1)}}\,.
\]
\item If $c_1 > 1, c_2 \leq 1,$ and $c_1 + c_2 \leq 2,$ then
\[
    \rho = \frac{1}{2}\,\frac{c_2(c_1 - c_2/2)}{\sqrt{(1 + c_2/2)p(c_1)}}\,.
\]
\item If $c_1 > 1, c_2 \leq 1,$ and $c_1 + c_2 > 2,$ then
\[
    \rho = -\frac{f(c_1, c_2)}{\sqrt{(1 + c_2/2)p(c_1)}}\,.
\]
\end{itemize}
In conclusion, the cases $c_2 = 2$ and $c_2 = 0,$ $c_1 > 0$ yield $\rho = 0$. 
\end{theorem}

\subsection{Asymptotic correlation for low regimes of \texorpdfstring{$d_1$}{d1}} \label{sec3.3} 
For the low regimes 
\[
    d_1 = o\left(n^{1/3}\log(n^{4/3})^{-7/3}\right),
\]
the replacement of $X_\inv^{(d_1)}$ with $\hat{X}_\inv^{(d_1)}$ is infeasible by Corollary \ref{cor3.5}, and we obtain Gaussian approximation for the original joint statistic $(X_\inv^{(d_1)}, X_\des^{(d_2)})^\top$. Therefore, we now investigate the correlation $\rho_n = \corr\left(X_\inv^{(d_1)}, X_\des^{(d_2)}\right)$ and its limit behavior. An interesting observation is that $X_\inv^{(d_1)}$ and $X_\des^{(d_1)}$ coincide if $d_1 = d_2 = 1,$ which implies $\rho = 1$ if $d_1 \equiv d_2 \equiv 1$ throughout. We will generally notice a nonzero limit correlation $\rho = \rho(d_1)$ if both $d_1, d_2$ remain fixed. If $d_1(n)$ is bounded but has multiple accumulation points, then $\rho$ does not exist. On the contrary, if $d_1(n)$ is not bounded, then we always have $\rho = 0$, in analogy with the sublinear cases in Theorem~\ref{thm3.10}. 

\subsection*{Symmetric groups} Similar as in Lemma~\ref{lemma3.10}, we provide an exact computation and classification of the covariance $\Cov(X_\inv^{(d_1)}, X_\des^{(d_2)})$ in the relevant regimes of $d_1$. 

\begin{theorem} \label{thm3.13}
For any two numbers $d_1, d_2 \in \{1, \ldots, n-1\}$ with $d_1 < n/2,$ the covariance of $X_\inv^{(d_1)}$ and $X_{\des}^{(d_2)}$ is given as follows, distinguishing by the relative order of $d_1$ and $d_2 \negmedspace:$ \par \bigskip \noindent 
If $d_1 < d_2,$ then $\Cov\left(X_\inv^{(d_1)}, X_{\des}^{(d_2)}\right) =  \Cov\left(\hat{X}_\inv^{(d_1)}, X_{\des}^{(d_2)}\right),$ i.e., 
\[
    \Cov\left(X_\inv^{(d_1)}, X_{\des}^{(d_2)}\right) = \begin{cases}
        d_1(d_1+1)/12, & d_1 + d_2 \leq n \\ (n(d_1 + d_2) - d_1d_2)/6 - (d_2(d_2 + 1) + n(n-1))/12, & d_1 + d_2 > n
    \end{cases}.
\]
If $d_1 = d_2 =: d < n/2,$ then
\[
    \Cov\left(X_\inv^{(d)}, X_{\des}^{(d)}\right) = \frac{n + d^2}{12}\,. 
\]
If $d_1 > d_2,$ then
\begin{align*}
    &\Cov\left(X_\inv^{(d_1)}, X_{\des}^{(d_2)}\right) = \frac{d_1d_2}{6} + \frac{n - d_2^2}{12}\,. 
\end{align*}
\end{theorem}

\begin{theorem} \label{thm3.14}
Let $d_1 = o\left(n^{1/3}\log(n^{4/3})^{-7/3}\right)$ and let $d_2 = d_2(n) \in \{1, \ldots, n-1\}$ be arbitrary. Then, for the limiting correlation 
\[
    \rho = \lim_{n\rightarrow \infty} \corr\left(X_\inv^{(d_1)}, X_\des^{(d_2)}\right),
\]
the following holds:
\begin{enumerate}
    \item[\emph{a)}] If $d_1 < d_2$ throughout, then $\rho = 0$.
    \item[\emph{b)}] If $d_1 \geq d_2$ throughout and $d_1$ eventually remains constant, then $\rho = 1/\sqrt{d_1}$. If $d_1$ is bounded but has multiple accumulation points, then $\rho$ does not exist.
    \item[\emph{b)}] If $d_1 \geq d_2$ and $d_1 \rightarrow \infty$ as $n \rightarrow \infty,$ then $\rho = 0$.
\end{enumerate}
\end{theorem}

In conclusion, if $d_1$ and $d_2$ both remain bounded, there is a sharp threshold between the cases $d_1 < d_2$ and $d_1 \geq d_2$. 

\subsection*{Signed and even-signed permutation groups} While the generalized inversions and descents statistics on the signed and even-signed permutation groups have a more complex structure, the asymptotic behavior of $\corr(X_\inv^{(d_1)}, X_\des^{(d_2)})$ is basically the same as for the symmetric groups. 

\begin{theorem} \label{thm3.15}
Consider the families of signed permutation groups $(B_n)_{n \in \N}$ or even-signed permutation groups $(D_n)_{n \in \N}$. Let $d_1 = o\left(n^{1/3}\log(n^{4/3})^{-7/3}\right)$ and let $d_2 = d_2(n) \in \{1, \ldots, d_{\max}-1\}$ be arbitrary. Then, for the limiting correlation 
\[
    \rho = \lim_{n\rightarrow \infty} \corr\left(X_\inv^{(d_1)}, X_\des^{(d_2)}\right),
\]
the same results as in Theorem \emph{\ref{thm3.14}} hold:
\begin{enumerate}
    \item[\emph{a)}] If $d_1 < d_2$ throughout, then $\rho = 0$.
    \item[\emph{b)}] If $d_1 \geq d_2$ throughout and $d_1$ eventually remains constant, then $\rho = 1/\sqrt{d_1}$. If $d_1$ is bounded but has multiple accumulation points, then $\rho$ does not exist.
    \item[\emph{b)}] If $d_1 \geq d_2$ and $d_1 \rightarrow \infty$ as $n \rightarrow \infty,$ then $\rho = 0$.
\end{enumerate}
\end{theorem}

\section{Extreme Value Limit Theorems}
\label{section4}

\subsection*{Extreme values of the univariate statistics \texorpdfstring{$X_\inv^{(d)}$}{Xinvd} and \texorpdfstring{$X_\des^{(d)}$}{Xdes}} 
We now postulate the univariate and bivariate extreme value limit theorems (EVLTs) for $X_\inv^{(d)}$ and $X_\des^{(d)}$ with the help of the methods used in \cite{dorr2025joint}. 
Recall (e.g., by \cite{leadbetter2012extremes}, Theorem 1.5.3) that by use of 
\begin{align*}
    \alpha_n &= \frac{1}{\sqrt{2\log n}}\,, & \beta_n &= \frac{1}{\alpha_n} - \frac{1}{2}\alpha_n\bigl(\log \log n + \log(4\pi)\bigr)\,,
\end{align*}
the maximum of i.i.d.\ variables $N_1, \ldots, N_n \sim \mathrm{N}(0,1)$ is attracted to the Gumbel distribution, that is, $\forall x \in \R$:
\[\PP\left(\max_{i=1, \ldots, n} N_i \leq \alpha_n x + \beta_n\right) \longrightarrow \exp(-\exp(-x)) =: \Lambda(x).\]
Since the random numbers of generalized inversions and descents are discrete distributions, we are interested in the row-wise extremes of a triangular array $(X_{n1}, \ldots, X_{nk_n}),$ where for each $n \in \N,$ the block $X_{n1}, \ldots, X_{nk_n}$ consists of i.i.d.\ samples based on one of the classical Weyl groups $S_n, B_n$ or $D_n$. Each of the following EVLTs states a regime for the number of samples $k_n$ for which Gumbel attraction is guaranteed.

For a univariate triangular array consisting of generalized descents, it is not necessary to use the H\'{a}jek projection. It has already been argued in \cite[Remark 4.2]{dorr2025joint} that a subexponential bound on $k_n$ can be obtained if the H\'{a}jek projection is not needed. So, in the univariate EVLT for generalized descents, we can proceed in analogy to \cite[Remark 4.2]{dorr2025joint}:

\begin{theorem} \label{thm4.1}
Let $(X_{nj})_{j=1,\ldots,k_n}$ be a row-wise i.i.d.\ triangular array with $X_{n1} \overset{\D}{=} X_\des^{(d)}$ for a sequence $d = d(n)$ with $1 \leq d \leq d_{\max}-1$ and $N_{n,d}^= \longrightarrow \infty,$ and let $M_n := \max\{X_{n1}, \ldots, X_{nk_n}\}$. Let $a_n := \sigma\bigl(X_\des^{(d)}\bigr)\alpha_{k_n}$ and let $b_n := \sigma\bigl(X_\des^{(d)}\bigr)\beta_{k_n} + \E\bigl(X_\inv^{(d)}\bigr).$ If $k_n = \exp\left(o\bigl((N_{n,d}^=)^{1/7}\bigr)\right),$ then
\[
    \forall x \in \R \negmedspace: \quad \PP(M_n \leq a_nx + b_n) \longrightarrow \Lambda(x)\,. 
\]
\end{theorem}

\begin{proof}
    By \eqref{4}, \eqref{2.7c}, \eqref{2.7d}, we can write each $X_{n1}$ as a sum of $N_{n,d}^=$ indicator variables with $\Delta_n \leq 3$ and $\Delta_n^* \leq 9$. This means we apply Theorem~\ref{thm3.3} with respect to the numbers of summands $N_{n,d}^=$ and for high-dimensional vectors built from $k_n$ iterations of $X_{n1},$ yielding
    \[
        |\PP(M_n \leq a_nx + b_n) - \PP(\mathcal{M}_n \leq \alpha_nx + \beta_n)| = O\left((N_{n,d}^=)^{-1/6}\log(k_n)^{7/6}\right),
    \]
    where $\mathcal{M}_n$ is the maximum of $n$ i.i.d.\ copies of the standard normal distribution. Thus we need to ensure
    \[
        (N_{n,d}^=)^{-1/6}\log(k_n)^{7/6} = o(1)\,,
    \]
    which exactly corresponds with the assumption of $k_n = \exp\left(o\bigl((N_{n,d}^=)^{1/7}\bigr)\right)$.
\end{proof}

A similar statement applies for $X_\inv^{(d)}$ if $d$ grows slowly enough. 

\begin{theorem} \label{thm3.6}
Let $(X_{n1}, \ldots, X_{nk_n})$ be a row-wise i.i.d. triangular array with $X_{n1} \overset{\D}{=} X_\inv^{(d)},$ where $d = d(n) = o(n^{1/3})$. Let $M_n, a_n, b_n$ be given in analogy to Theorem \emph{\ref{thm4.1}}. If $k_n = \exp\bigl(o(n^{1/7}d^{-3/7})\bigr),$ then
\[
    \forall x \in \R \negmedspace: \quad \PP(M_n \leq a_nx + b_n) \longrightarrow \Lambda(x)\,. 
\]
\end{theorem}

\begin{proof}
According to the above considerations, the maximum degrees $\Delta_n, \Delta_n^*$ in the dependency graphs of the representations \eqref{3}, \eqref{2.7a}, and \eqref{2.7b} are bounded in the way of $\Delta_n \leq 4d$ and $\Delta_n^* \leq 8d$. These representations are based on $\Theta(nd)$ summands. By analogy with the proof of Theorem~\ref{thm4.1}, an application of Theorem~\ref{thm3.3} with $k_n$ i.i.d.\ iterations of $X_\inv^{(d)}$ yields  
\[
    |\PP(M_n \leq a_nx + b_n) - \PP(\mathcal{M}_n \leq \alpha_nx + \beta_n)| = O(n^{-1/6}d^{1/2})\log(k_n)^{7/6}\,.
\]
Again, plugging $k_n = \exp\bigl(o(n^{1/7}d^{-3/7})\bigr)$ according to the assumption yields
\[
    |\PP(M_n \leq a_nx + b_n) - \PP(\mathcal{M}_n \leq \alpha_nx + \beta_n)| = o(1)\,,
\]
from which the claim follows.
\end{proof}

For any other growth rate of $d,$ we can state an EVLT only for the cases covered by Lemma \ref{lemma3.4}, and due to the use of H\'{a}jek's projection, we can only impose a strongly reduced asymptotic bound on $k_n$.

\begin{theorem} \label{thm4.3}
Let $(X_{n1}, \ldots, X_{nk_n})$ be a row-wise i.i.d. triangular array with $X_{n1} \overset{\D}{=} X_\inv^{(d)}$ and with $d = d(n)$ such that $d = \omega(n^{1/2})$. Assume $k_n \log(k_n) = o(d^2/n)$ and let $M_n, a_n, b_n$ be given in analogy to Theorem \emph{\ref{thm4.1}}. Then,
\[
    \forall x \in \R \negthickspace: \quad \PP(M_n \leq a_nx + b_n) \longrightarrow \Lambda(x)\,.
\]
\end{theorem}

\begin{proof}
The conditions on $d$ ensure that $1 - \Var(X_\inv^{(d)})/\Var(\hat{X}_\inv^{(d)}) = o(1),$ but the rate of convergence determines the bound on $k_n$ by means of \cite[Eq. (9), (10)]{dorr2025joint}. We compute this rate for symmetric groups, since the same conclusions can be obtained on the other classical Weyl groups. By the proof of Lemma~\ref{lemma3.4}, we have for $d \leq n/2$ that
\begin{align*}
    1 - \frac{\Var(X_\inv^{(d)})}{\Var(\hat{X}_\inv^{(d)})} &= 1 - \frac{4d^3 + 6nd + 3d^2 - d}{4d^3 + 6d^2 + 2d} = \frac{4d^2 + 6n + 3d - 1}{4d^2 + 6d + 2} \\
    &= 1 - \frac{4d^2}{4d^2 + 6d + 2} - \frac{3d+1}{4d^2 + 6d + 2} - \frac{6n}{4d^2 + 6d + 2} \\
    &= \Theta\left(\frac{1}{d}\right) - \frac{6n}{4d^2 + 6d + 2}.
\end{align*}
Apparently, $6n/(4d^2 + 6d + 2)$ always dominates $1/d$. In conclusion, 
\[
    1 - \frac{\Var(X_\inv^{(d)})}{\Var(\hat{X}_\inv^{(d)})} = \Theta\left(\frac{n}{d^2}\right),
\]
giving the condition of $k_n \log(k_n) = o(d^2/n)$ according to the arguments in the proof of \cite[Theorem 4.1]{dorr2025joint}. From here, we proceed as in the proof of \cite[Theorem 4.1]{dorr2025joint}. In the case of $d > n/2,$ we always have $\Var(X_\inv^{(d)}) \sim \Var(\hat{X}_\inv^{(d)})$ and therefore, 
\[1 - \frac{\Var(X_\inv^{(d)})}{\Var(\hat{X}_\inv^{(d)})} = \Theta\left(\frac{1}{d}\right) = \Theta\left(\frac{1}{n}\right) = \Theta\left(\frac{n}{d^2}\right).\]
Again, the proof now follows the same steps as in \cite[Theorem 4.1]{dorr2025joint}.
\end{proof} 

\subsection*{Extreme values of the joint distribution of \texorpdfstring{$X_\inv^{(d)}$}{Xinvd} and \texorpdfstring{$X_\des^{(d)}$}{Xdesd}} 

We recall (e.g., by \cite[Theorem 3]{sibuya1960bivariate}) that any bivariate normal distribution with correlation $\rho < 1$ is max-attracted to the two-dimensional Gumbel distribution with independent marginals. Let $\mathbf{N}_1, \ldots, \mathbf{N}_n$ be i.i.d.\ samples of $N_2\left(0, \begin{pmatrix} 1 & \rho \\ \rho & 1 \end{pmatrix} \right)$. Taking $\boldsymbol{\upalpha}_n := (\alpha_n, \alpha_n),$ $\boldsymbol{\upbeta}_n := (\beta_n, \beta_n) \in \R^2$ and writing ''$*$'' for component-wise multiplication, we have by \cite[Theorem 3]{sibuya1960bivariate} that for $\mathbf{x} = (x_1, x_2) \in \R^2$:
\[
    \PP\left(\max_{i=1, \ldots, n} \mathbf{N}_i \leq \boldsymbol{\upalpha}_n * \mathbf{x} + \boldsymbol{\upbeta}_n\right) \longrightarrow \exp\bigl(-\exp(-x_1) - \exp(-x_2)\bigr) =: \Lambda_2(x).
\]
Now, the bivariate EVLT for $\left(X_\inv^{(d_1)}, X_\des^{(d_2)}\right)^\top$ can be stated analogously to the univariate EVLTs. Note that the second component $X_\des^{(d_2)}$ does not interfere with the arguments in the proofs of Theorems ~\ref{thm3.6} and~\ref{thm4.3}. However, we assume that the two sequences $d_1 = d_1(n),$ $d_2 = d_2(n)$ with $d_1, d_2 \leq d_{\max} - 1$ $\forall n \in \N$ are chosen in such a way that the limiting correlation $\rho = \lim_{n\rightarrow \infty} \corr\left(X_\inv^{(d_1)}, X_\des^{(d_2)}\right)$ exists, see Theorems~\ref{thm3.10}, \ref{thm3.12}, \ref{thm3.14}, \ref{thm3.15}. 

\begin{remark}
The only scenario for which $\rho = 1$ occurs if $d_1, d_2 \equiv 1$ throughout, i.e., both $X_\inv^{(d_1)}$ and $X_\des^{(d_2)}$ are the common descent statistic $X_\des$. In this case, the bivariate EVLT trivially reduces to the established EVLT for $X_\des,$ see \cite[Theorem 13]{dorr2022extreme}. Hence, we ignore this scenario.
\end{remark}

\begin{theorem} \label{thm4.2}
Let $d_1 = d_1(n),$ $d_2 = d_2(n)$ be two sequences as above, and let $(X_{n1}, \ldots, X_{nk_n})$ be a row-wise i.i.d.\ triangular array with $X_{n1} \overset{\D}{=} \left(X_\inv^{(d_1)}, X_\des^{(d_2)}\right)^\top$. Then: 
\begin{enumerate}
\item[\emph{(a)}] If $d_1 = o(n^{1/3}),$ we assume $k_n = \exp\left(o\left((N_{n,d_2}^=)^{1/7} \wedge n^{1/7}d_1^{-3/7}\right)\right).$
\item[\emph{(b)}] If $d_1 = \omega(n^{1/2}),$ we assume $k_n = o(d^2/n)$. 
\end{enumerate}
Let $M_n := \max\{X_{n1}, \ldots, X_{nk_n}\}$ be the row-wise maximum. Let $\mu_n := \E(X_{n1})$ and let $s_n := \bigl(\sigma(X_\inv^{(d)}), \sigma(X_\des)\bigr)$. Let $a_n := s_n * \boldsymbol{\upalpha}_{k_n}$ and let $b_n := s_n * \boldsymbol{\upbeta}_{k_n} + \mu_n$. Then,
\[\forall \mathbf{x} \in \R^2 \negthickspace: \quad \PP(M_n \leq a_n * \mathbf{x} + b_n) \longrightarrow \Lambda_2(\mathbf{x})\,.\]
\end{theorem}

\section{Conclusion and Outlook} \label{sec5}

We successfully investigated the bivariate CLT for the generalized inversion and descent statistics $X_\inv^{(d_1)}$ and $X_\des^{(d_2)}$ on all classical Weyl groups. For low regimes $d_1 \approx o(n^{1/3})$, the CLT was achieved through direct Gaussian approximation, while for larger regimes $d_1 \gg n^{1/2}$, the CLT was achieved with help of the H\'{a}jek projection $\hat{X}_\inv^{(d_1)}$. Corresponding extreme value results for these two statistics were derived as well. There is a remaining gap consisting of the regimes $n^{1/3} \ll d_1 \ll n^{1/2}$ which could be closed upon availability of suitable Gaussian approximation results. Simulations suggested that the CLT applies for these regimes as well.

The investigation of the asymptotic correlation $\rho$ showed several interesting phenomena. When $d_1$ diverges to infinity as $n\rightarrow \infty$ but remains sublinear, the correlation is always zero, with the exception of some ''complementary'' cases where $d_2$ is chosen so that $d_1 + d_2 = n(1 + o(1))$. If $d_2$ is sublinear, $\rho = 0$ applies as well. If both $d_1$ and $d_2$ are linear, then $\rho$ is an explicit function of the limits of $d_1/n$ and $d_2/n$. On the other hand, if $d_1$ remains bounded, there is a clear distinction between the cases $d_1 < d_2$ and $d_1 \geq d_2$. The former case guarantees $\rho = 0,$ while in the latter, the CLT can only be satisfied if $d_1$ eventually remains constant, in which case $\rho = 1/\sqrt{d_1} > 0$. 

There are several related questions beyond the scope of this paper. Other well-known permutation statistics can be generalized in a similar pattern. For example, on the symmetric groups, one can consider the \textit{generalized major index}
\[
    X_{\mathrm{maj}}^{(d)} := \sum_{i=1}^{n-d} i \cdot \textbf{1}\{Z_i > Z_{i+d}\}
\]
and the \textit{generalized peaks}
\[
    X_{\mathrm{peak}}^{(d,e)} := \sum_{i=d+1}^{n-e} \textbf{1}\{Z_{i-d} < Z_i, Z_{i+e} < Z_i\}\,.
\]
While the common major index $X_\mathrm{maj}$ is equidistributed to the common inversion statistic $X_\inv,$ this no longer applies for the generalized counterparts. Further generalized statistics arise by considering joint features of a permutation and its inverse. The most basic example is the \textit{two-sided Eulerian statistic} $X_T(\pi) := X_\des(\pi) + X_\des(\pi^{-1})$ for $\pi \in S_n$ (or $B_n$, $D_n$, see \cite{bruck2019central, chatterjee2017central, rottger2018asymptotics}). Its generalized version can be defined as $X_T^{(d, e)}(\pi) := X_\des^{(d)}(\pi) + X_\des^{(e)}(\pi^{-1})$. Even for the original statistic $X_T$, the joint asymptotic normality with $X_\inv$ is not investigated yet due to the complex dependency structure of $X_T$. 

Besides further generalized statistics, one can also drop the basic framework of drawing elements of classical Weyl groups uniformly at random. Prominent examples of non-uniform distributions on classical Weyl groups are the Ewens distribution (cf.\  \cite{feray2012asymptotic}) and the Mallows distribution (cf.\ \cite{he2022central, sun2024central}). 

\begin{appendices}

\section{Calculations of H\'{a}jek projections} \label{appendixB}

\subsection*{Proof of Lemma~\ref{lemma3.4}}
\begin{figure}[h]
    \centering
    \includegraphics[width=\textwidth]{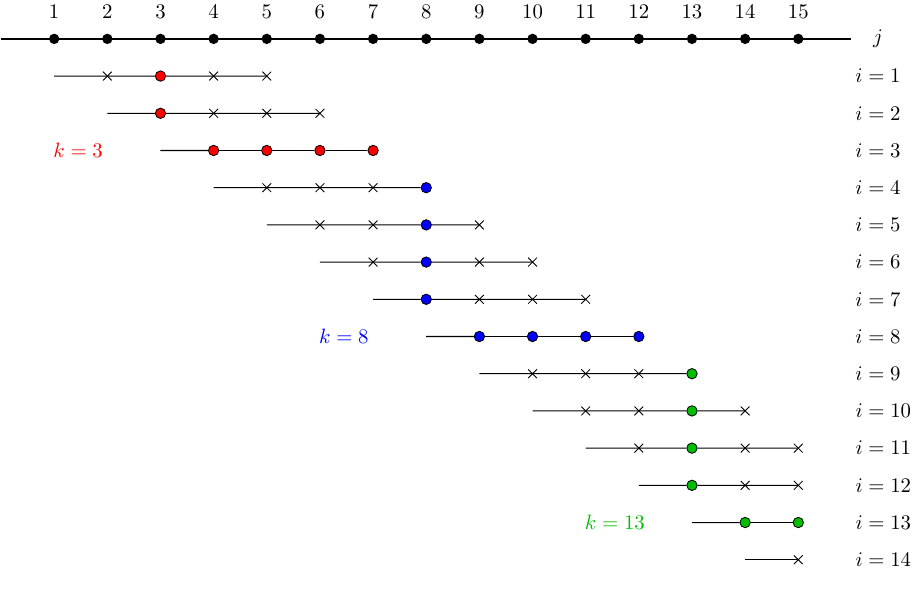}
    \caption{Overview of relevant pairs $(i,j)$ for computing the variance of $\E\left(X_\inv^{(d)} \mid Z_k\right)$, where $n=15$ and $d=4$. For each of the regions $K_1, K_2, K_3,$ an exemplary index $k$ is chosen, and the pairs that give a non-trivial contribution are highlighted in red for $k \in K_1$, in blue for $k \in K_2$, and in green for $k \in K_3$.}
    \label{fig5.2}
\end{figure}

Consider $\hat{X}_\inv^{(d)}$ on symmetric groups. By \eqref{3}, we have
\[
    \E(X_\inv^{(d)} \mid Z_k) = \sum_{(i,j) \in \NN_{n,d}} \PP(Z_i > Z_j \mid Z_k)\,,  ~~~\text{where}\,\,\,\PP(Z_i > Z_j \mid Z_k) = \begin{cases} 1/2, & k \notin \{i,j\} \\ Z_k, & k=i \\ 1 - Z_k, & k=j \end{cases}.
\]
Only the pairs $(i,j)$ with $k \in \{i,j\}$ contribute to $\Var\left(\E\left(X_\inv^{(d)} \mid Z_k\right)\right)$. The number of these \textit{non-trivial pairs} depends on whether $k$ belongs to $K_1, K_2,$ or $K_3$, see Remark~\ref{bem2.2}. Figure~\ref{fig5.2} visualizes this case distinction for the exemplary choice of $n=15$ and $d=4$.

\begin{figure}[h]
    \centering
    \includegraphics[scale=0.8]{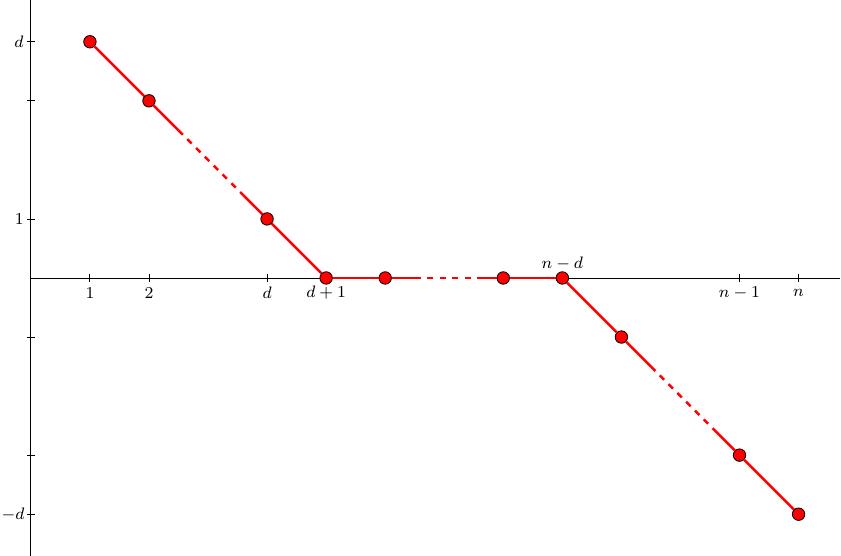}
    \caption{Plot of $\omega_d(k)$ in the case of $d \leq n/2$.}
    \label{fig4.6}
\end{figure}

\begin{figure}[h]
    \centering
    \includegraphics[scale=0.7]{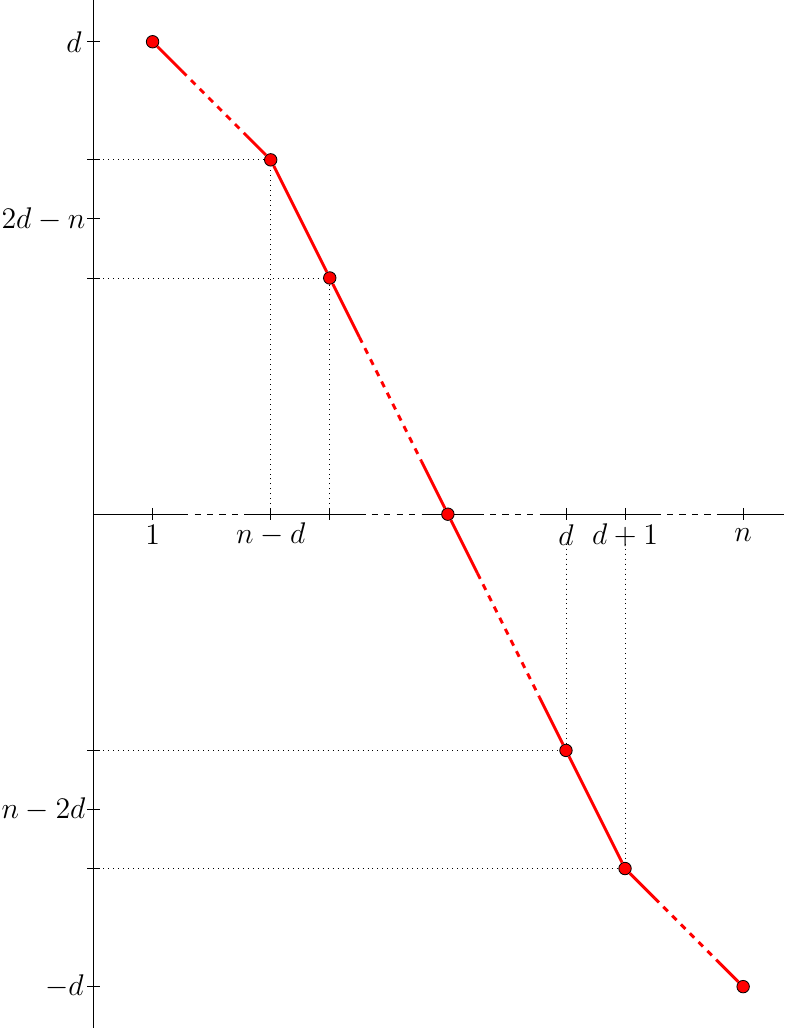}
    \caption{Plot of $\omega_d(k)$ in the case of $d > n/2$.}
    \label{fig4.7}
\end{figure}

\smallskip
\noindent We first assume $n \leq d/2$. If $k \in K_2,$ then the non-trivial parts are
\[
    \underbrace{Z_k + Z_k + \ldots + Z_k}_{d~\text{times}} + \underbrace{(1 - Z_k) + (1 - Z_k) + \ldots + (1 - Z_k)}_{d~\text{times}} = d.
\]
This means $\E\left(X_\inv^{(d)} \mid Z_k\right)$ is constant due to cancellation, and vanishes when computing the variance. Otherwise, we find
\[
    \E\left(X_\inv^{(d)} \mid Z_k\right) = \begin{cases}
        (k-1)(1 - Z_k) + dZ_k + \text{const} = (d+1-k)Z_k + \text{const}, & k \in K_1 \\ (n-k)Z_k + d(1 - Z_k) + \text{const} = (n-d-k)Z_k + \text{const}, & k \in K_3
    \end{cases}.
\]
So, the overall representation of $\hat{X}_\inv^{(d)}$ for the case $d \leq n/2$ is
\begin{subequations}
\begin{equation}
    \hat{X}_\inv^{(d)} = \sum_{k=1}^n \omega_d(k)Z_k + \text{const}\,, \quad \text{with} \quad \omega_d(k) := \begin{cases} d-k+1\,, & k \in K_1 \\ 0, & k \in K_2 \\ n-d-k\,, & k \in K_3 \end{cases}. \label{omegashort}
\end{equation}
Therefore,
\begin{align*}
    \Var\left(\hat{X}_\inv^{(d)}\right) &= \Var\left(\sum_{k=n-d+1}^n (n-k-d)Z_k + \sum_{k=1}^d (d+1-k)Z_k\right) \\
    &= \sum_{k=n-d+1}^n \frac{1}{12}(n-k-d)^2 + \sum_{k=1}^d \frac{1}{12}(d+1-k)^2 \\
    &= \frac{1}{12}\sum_{k=1}^d k^2 + \frac{1}{12}\sum_{k=1}^d (-k)^2 = \frac{1}{6}\frac{d(d+1)(2d+1)}{6} \\
    &= \frac{1}{72}\left(4d^3 + 6d^2 + 2d\right). 
\end{align*}
So, the leading term is always $4d^3/72$. In light of Theorem~\ref{thm2.3}, according to which $\Var\left(X_\inv^{(d)}\right)$ contains the monomials $6nd/72$ and $4d^3/72,$ the application of Theorem~\ref{thm2.8} requires the condition  $d^3 \gg nd \Longleftrightarrow d \gg \sqrt{n}$, as given by the claim. \par 
\bigskip
Now, we consider the case of $d > n/2$. We then have $n-d < d$, and the regions $K_1, K_2, K_3$ are redefined according to Remark~\ref{bem2.2}. For the non-trivial parts, we note that:
\begin{itemize}
    \item If $k \in K_1$ or $k \in K_3,$ then the non-trivial parts are the same as above.
    \item If $k \in K_2,$ then the non-trivial parts yield $\E\left(X_\inv^{(d)} \mid Z_k\right) = (n-2k+1)Z_k + \text{const}$.
\end{itemize}
So, we again obtain a representation in the way of
\begin{equation}
    \hat{X}_\inv^{(d)} = \sum_{k=1}^n \omega_d(k)Z_k + \text{const}\,, \qquad \text{with} \quad \omega_d(k) := \begin{cases} d-k+1, & k \in K_1 \\ n-2k+1, & k \in K_2 \\ n-d-k, & k \in K_3 \end{cases}. \label{omegalong} 
\end{equation}
\end{subequations}
For illustration purposes, the coefficients $\omega_d(k)$ appearing in \eqref{omegashort} and \eqref{omegalong} are displayed in Figures \ref{fig4.6} and~\ref{fig4.7}, respectively.

From \eqref{omegalong}, we state that
\begin{subequations}
\begin{align}
    \Var\left(\hat{X}_\inv^{(d)}\right) &= \frac{1}{12} \sum_{k=1}^{n-d} (d+1-k)^2 \label{6a} \\
    &+ \frac{1}{12} \sum_{k=n-d+1}^d (n-2k+1)^2 \label{6b} \\
    &+ \frac{1}{12}\sum_{k=d+1}^n (n-k-d)^2 \label{6c} . 
\end{align}
\end{subequations}
By appropriate index shifting, we calculate 
\begin{align*}
    \eqref{6a}, \eqref{6c} &= \frac{1}{12}\left(\sum_{k=1}^d k^2 - \sum_{k=1}^{2d-n} k^2\right) \\
    \Longrightarrow \eqref{6a} + \eqref{6c} &= \frac{1}{6}\left(\sum_{k=1}^d k^2 - \sum_{k=1}^{2d-n} k^2\right) \\
    &= \frac{1}{36}(n-d)\bigl(14d^2 + d(9 - 10n) + 2n^2 - 3n + 1\bigr), \\
    \eqref{6b} &= \frac{1}{36}(2d-n)(4d^2 - 4dn + n^2-1). 
\end{align*}
This gives the total result
\begin{align*}
    \Var\left(\hat{X}_\inv^{(d)}\right) = \frac{1}{36}\left(-6d^3+12d^2n - 9d^2 - 6dn^2 + 12dn - 3d + n^3 - 3n^2 + 2n\right).
\end{align*}
In contrast, by Theorem~\ref{thm2.3}, \[\Var\left(X_\inv^{(d)}\right) = \frac{1}{36}\left(-6d^3 + \left(12n - \frac{21}{2}\right)d^2 - \left(6n^2 - 15n + \frac{9}{2}\right)d + n^3 - 3n^2 + 2n\right). \]
Since $n/2 < d < n$, all monomials of order 3 are leading terms. It is easily seen that these leading terms are matching, i.e.,
\[
\Var\left(X_\inv^{(d)}\right), \Var\left(\hat{X}_\inv^{(d)}\right) = \frac{1}{36}\left(-6d^3 + 12nd^2 - 6n^2d + n^3 + O(n^2)\right).
\]
This proves the lemma for symmetric groups. \qed

\subsection*{Proof of Lemma~\ref{lemma3.7}}

It is sufficient to give the proof only for the signed permutation groups $B_n,$ since the difference between $X_\inv^{B,(d)}$ and $X_\inv^{D,(d)}$ is asymptotically negligible (cf. \eqref{2.7a} and \eqref{2.7b}). 

\begin{figure}[ht]
    \centering
    \includegraphics[scale=1]{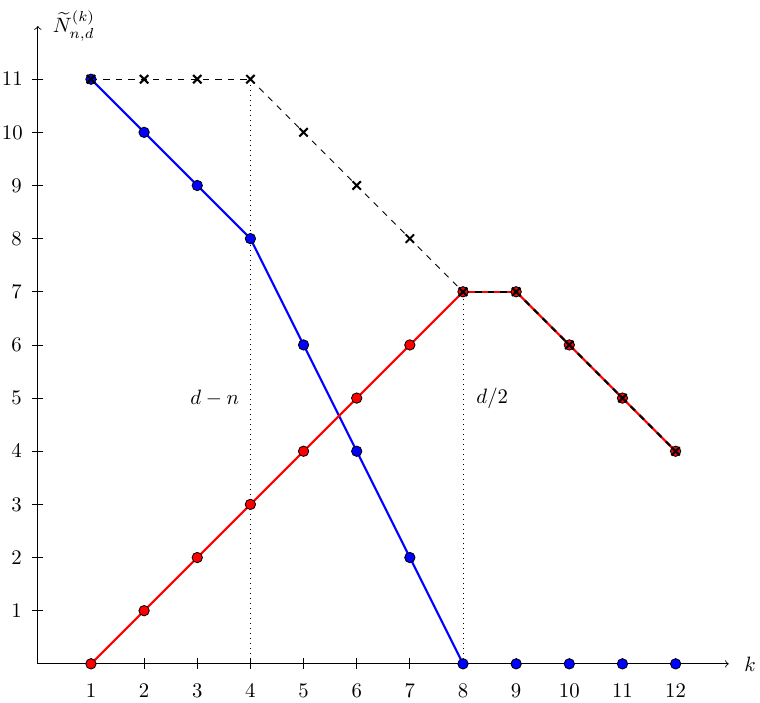}
    \caption{Plots of the numbers of pairs $(i,k)$ (red) and $(k,j)$ (blue) in $\widetilde{\NN}_{n,d}$. The sum of these two numbers is $\widetilde{N}_{n,d}^{(k)},$ which is displayed by the black crosses.}
    \label{fig4.9}
\end{figure}

To compute $\Var(\hat{X}_\inv^{B,(d)}),$ we ignore all constant parts appearing in $\hat{X}_\inv^{(d)}$. By \eqref{2.7a}, we have
\begin{align*}
\E\left(X_\inv^{B,(d)} \mid Z_k\right) &=
\sum_{(i,j) \in \NN_{n,d}} \PP(Z_i > Z_j \mid Z_k) + \sum_{(i,j) \in \widetilde{\NN}_{n,d}} \PP(-Z_i > Z_j \mid Z_k) \\
&\hspace{1.15cm}+ \sum_{i=1}^{n \wedge d} \PP(Z_i < 0 \mid Z_k)\,.
\end{align*}
Similar to the symmetric groups, we will compute coefficients $\omega_d(k)$ such that
\[
\sum_{(i,j) \in \NN_{n,d}} \PP(Z_i > Z_j \mid Z_k) + \sum_{(i,j) \in \widetilde{\NN}_{n,d}} \PP(-Z_i > Z_j \mid Z_k) = \omega_d(k)Z_k + \text{const}\,.
\]
If $d < k,$ then the third sum $\displaystyle{\sum\nolimits_{i=1}^{n \wedge d} \PP(Z_i < 0 \mid Z_k) = d/2}$ is constant. Otherwise, we have
\begin{align}
    \Var\left(\E\left(X_\inv^{B,(d)} \mid Z_k\right)\right) &= \Var\Bigl(\omega_d(k)Z_k + \textbf{1}\{Z_k < 0\} + \text{const}\Bigr) \nonumber \\
    &= \Var(\omega_d(k)Z_k) + \Var(\textbf{1}\{Z_k < 0\}) + 2\Cov(\omega_d(k)Z_k, \textbf{1}\{Z_k < 0\}) \nonumber \\
    &= \frac{\omega_d(k)^2}{3} + \frac{1}{4} + 2\omega_d(k)\underbrace{\Cov(Z_k, \textbf{1}\{Z_k < 0\})}_{=~-1/4}\,, \label{5.7}
\end{align} 
from which we see that even if $d < k,$ the leading terms of $\Var\left(\hat{X}_\inv^{B,(d)}\right)$ are not influenced by the third sum. In light of \eqref{5.7}, we have to determine the linear coefficients $\omega_d(k)$ stemming from 
\[
\sum_{(i,j) \in \NN_{n,d}} \PP(Z_i > Z_j \mid Z_k) + \sum_{(i,j) \in \widetilde{\NN}_{n,d}} \PP(-Z_i > Z_j \mid Z_k)\,.
\]
We write $\omega_d(k) = \omega_d(k)^+ + \omega_d(k)^-$, with $\omega_d(k)^+$  stemming from the first sum and $\omega_d(k)^-$ stemming from the second.
For $\displaystyle{\omega_d(k)^+},$ we can use the counting method previously used in the proof of Lemma~\ref{lemma3.4}. However, we have to take into account that, as $Z_i \sim U(-1,1)$,
\begin{align*}
    \PP(Z_k > Z_j \mid Z_k) &= \frac{Z_k+1}{2}\,, & \PP(Z_i > Z_k \mid Z_k) &= \frac{1-Z_k}{2}\,, \\[1ex]
    \PP(-Z_k > Z_j \mid Z_k) &= \frac{1-Z_k}{2}\,, & \PP(Z_i > Z_k \mid Z_k) &= \frac{1-Z_k}{2}\,.
\end{align*}
In conclusion, the coefficients $\displaystyle{\omega_d(k)^+}$ on $B_n$ are half of the coefficients $\omega_d(k)$ on $S_n$ if $d < n$. Otherwise, for $d \geq n$, we always have $\omega_d(k)^+ = (n-2k+1)/2$. Moreover, $\omega_d(k)^- = -\widetilde{N}_{n,d}^{(k)}/2,$ where
\[
    \widetilde{N}_{n,d}^{(k)} := |\{(i,j) \in \widetilde{\NN}_{n,d} \mid \text{either}~ i=k~\text{or}~j=k\}|\,.
\]
For all $d \in \{1, \ldots, d_{\max}\},$ we can find that
\[
    \widetilde{N}_{n,d}^{(k)} = \begin{cases} n - 1, & k \leq d-n \\ d - k - 1, & d - n < k \leq d/2 \\ d - k, & d/2 < k \leq d \\ 0, & k > d\end{cases}\,.
\]
The first case $k \leq d-n$ is redundant if $d \leq n,$ and the last case $k > d$ is redundant if $n \geq d$. We give illustrations for the computation of $\widetilde{N}_{n,d}^{(k)}$ in the case $d > n$. Figure~\ref{fig4.9} displays the number of pairs $(i,k)$ and the number of pairs $(k,j)$ for the exemplary case of $n = 12, d = 16$. 

By analogy with the proofs of \cite[Theorems A.4 and A.13]{meier2022central}, we need to distinguish the four cases $d \leq n/2,$ $n/2 \leq d \leq 2n/3$, $2n/3 \leq d < n$, and $d \geq n$. If $d \leq n/2$, then all pairs in $\widetilde{\NN}_{n,d}$ are located within $K_1$, yielding
\begin{subequations}
\begin{align}
    \omega_d(k) &= \begin{cases} (d-k+1)/2 - (d-k-1)/2, & k \leq d/2 \\ (d-k+1)/2 - (d-k)/2, & d/2 < k \leq d \\ 0, & d < k \leq n-d \\ (n-d-k)/2, & n-d < k \leq n \end{cases} \nonumber \\[1ex] 
    &= \begin{cases} 1, & k \leq d/2 \\ 1/2, & d/2 < k \leq d \\ 0, & d < k \leq n-d \\ (n-d-k)/2, & k > n-d \end{cases} \label{omega_a} \\[1ex]
    &= \begin{cases} (n-d-k)/2, & k > n-d \\ 0, & d < k \leq n - d \\ O(1), & \text{otherwise} \end{cases}\,. \nonumber
\end{align}
In conclusion, if $d \leq n/2,$ then
\begin{align*}
    \Var\left(\hat{X}_\inv^{B,(d)}\right) &= \frac{1}{3}\sum_{k=1}^n \omega_d(k)^2 + O(d^2) = \frac{1}{12}\sum_{k=n-d+1}^n (n-d-k)^2 + O(d^2) \\
    &= \frac{1}{12}\frac{d(d+1)(2d+1)}{6} + O(d^2) = \frac{1}{36}d^3 + O(d^2).
\end{align*}
Due to $\Var\left(X_\inv^{B,(d)}\right) = d^3/36 + nd/12 + O(d^2)$ according to Lemma~\ref{lemma2.3.8}, we again obtain the  condition $d \gg \sqrt{n}$. \par 
\bigskip
\noindent If $n/2 < d < n$, then the pairs in $\widetilde{\NN}_{n,d}$ also cover $K_2$. For $k \in K_2,$ there cannot be any pairs $(k,j)$ if $n/2 < d \leq 2n/3,$ while this is possible if $d > 2n/3$. However, the difference between these two subcases is only marginal. If $n/2 < d \leq 2n/3,$ we obtain
\begin{equation}
    \omega_d(k) = \begin{cases} 1, & k \leq d/2 \\ 1/2, & d/2 < k \leq n - d \\ (n-d-k+1)/2, & n - d < k \leq d \\ (n-d-k)/2, & k > d\end{cases}\,. \label{omega_b}
\end{equation}
If $d > 2n/3,$ then
\begin{equation}
    \omega_d(k) = \begin{cases} 1, & k \leq n-d \\ (n-d-k+2)/2, & n-d < k \leq d/2 \\ (n-d-k+1)/2, & d/2 < k \leq d \\ (n-d-k)/2, & k > d\end{cases}\,. \label{omega_c}
\end{equation}
In conclusion, if $n/2 < d < n,$ then we obtain, as before,
\begin{align*}
    \Var\left(\hat{X}_\inv^{B,(d)}\right) &= \frac{1}{3}\sum_{k=1}^n \omega_d(k)^2 + O(d^2) \\
    &= \frac{1}{12}\sum_{k=n-d+1}^d (n-d-k+1)^2 + \frac{1}{12}\sum_{k=d+1}^n (n-d-k)^2 + O(d^2) \\
    &= \frac{1}{12}\sum_{k=n-d+1}^d (n-d-k+1)^2 + \frac{1}{12}\sum_{k=1}^{2d-n} k^2 + O(d^2) \\
    &= \frac{1}{12}\sum_{k=1}^{2d-n-1} k^2 + \frac{1}{12}\left(\sum_{k=1}^d k^2 - \sum_{k=1}^{2d-n} k^2\right) + O(d^2) \\
    &= \frac{1}{36}d^3 + O(d^2)\,. \qedhere
\end{align*}

    \label{fig4.8}

\noindent At last, if $d \geq n$, we have
\begin{equation}
    \omega_d(k) = \begin{cases} 1-k, & 1 \leq k \leq d-n \\(n-d-k+2)/2, & d-n < k \leq d/2, \\ (n-d-k+1)/2, & d/2 < k \leq n. \end{cases} \label{omega_d}
\end{equation}
\end{subequations}
We compute
\begin{align*}
    \sum_{k=1}^n \omega_d(k)^2 &= \sum_{k=1}^{d-n} (k-1)^2 + \frac{1}{4}\sum_{k=d-n+1}^n (n-d+1-k)^2 + O(n^2) \\
    &= \frac{1}{6}\bigl(2d^3 - 3d^2(2n+1) + d(6n^2 + 6n + 1) - n(2n^2 + 3n + 1)\bigr) \\
    &\qquad + \frac{1}{24}(2n-d)(14d^2 - 20dn - 9d + 8n^2 + 6n + 1) + O(n^2)\,.
\end{align*}
Due to $\Var(Z_k) = 1/3,$ we obtain
\begin{align*}
    \Var(\hat{X}_\inv^{B,(d)}) &= \frac{1}{72}\Bigl(-6d^3 + 3d^2(8n-1) - 3d(8n^2 - 1) + 2n(4n^2 - 1)\Bigr) + O(n^2) \\
    &= -\frac{1}{12}d^3 + \frac{1}{3}nd^2 - \frac{1}{3}n^2d + \frac{1}{9}n^3 + O(n^2)\,.
\end{align*}
By Lemma~\ref{lemma2.3.8}, this also applies for $\Var(X_\inv^{B,(d)})$, completing the proof. \qed

\section{Calculations of covariances and correlations involving H\'{a}jek projections} \label{appendixC}

\subsection*{Proof of Lemma \ref{lemma3.10}}

Recall that by \eqref{omegashort} and \eqref{omegalong}, we have
\[
    \hat{X}_\inv^{(d_1)} = \sum_{k=1}^d \omega_{d_1}(k)Z_k
\]
with $Z_1, \ldots, Z_n \overset{\text{i.i.d.}}{\sim} U(0,1)$ and 
\begin{align*}
    d_1 \leq n/2 &\Longrightarrow \omega_{d_1}(k) = \begin{cases} d_1 - k + 1\,, & k \leq d_1 \\ 0, & d_1 < k \leq n - d_1 \\ n - d_1 - k\,, & k > n - d_1 \end{cases}, \\
    d_1 > n/2 &\Longrightarrow \omega_{d_1}(k) = \begin{cases} d_1 - k +1, & k \leq n - d_1 \\ n - 2k + 1, & n - d_1 < k \leq d_1 \\ n- d_1 - k, & k > d_1 \end{cases}.
\end{align*}
Therefore, 
\[
    \Cov(\hat{X}_\inv^{(d_1)}, X_\des^{(d_2)}) = \sum_{j=1}^n \omega_{d_1}(j)\,\Cov\left(Z_j, X_\des^{(d_2)}\right).
\]
Consider a fixed $j \in \{1, \ldots, n\}$. For independence reasons, we have with slight abuse of notation:
\[
    \Cov\left(Z_j, X_\des^{(d_2)}\right) = \Cov\left(Z_j, \textbf{1}\{Z_{j-d_2} > Z_j\} + \textbf{1}\{Z_j > Z_{j+d_2}\}\right).
\]
By arguments in the proof of \cite[Lemma 2.6b)]{dorr2025joint}, we have
\begin{align*}
    \E(Z_j\textbf{1}\{Z_j < Z_{j-d_2}\}) = \frac{1}{6} &\Longrightarrow \Cov(Z_j, \textbf{1}\{Z_j < Z_{j-d_2}\}) = -1/12\,, \\
    \E(Z_j\textbf{1}\{Z_j > Z_{j+d_2}\}) = \frac{1}{3} &\Longrightarrow \Cov(Z_j, \textbf{1}\{Z_j > Z_{j+d_2}\}) = 1/12\,.
\end{align*}
If both indicators $\textbf{1}\{Z_{j-d_2} > Z_j\}$ and $\textbf{1}\{Z_j > Z_{j+d_2}\}$ are applicable, i.e., if both $j > d_2$ and $j \leq n - d_2,$ then $\Cov(Z_j, X_{\des}^{(d_2)}) = 0$. On the contrary, if $d_2 > n/2,$ then no indicator is applicable for all $j$ with $n-d_2 < j \leq d_2,$ also giving $\Cov(Z_j, X_{\des}^{(d_2)}) = 0$. In conclusion, we can write 
\[
    \Cov\left(\hat{X}_\inv^{(d_1)}, X_\des^{(d_2)}\right) = \sum_{j=1}^n \omega_{d_1}(j)c_{d_2}(j)\,,
\]
where
\begin{align*}
    d_2 \leq n/2 &\Longrightarrow c_{d_2}(j) = \begin{cases} +1/12\,, & j \leq d_2 \\ 0, & d_2 < j \leq n - d_2 \\ -1/12\,, & j > n-d_2 \end{cases}, \\
    d_2 > n/2 &\Longrightarrow c_{d_2}(j) = \begin{cases} +1/12\,, & j \leq n - d_2 \\ 0, & n - d_2 < j \leq d_2 \\ -1/12\,, & j > d_2 \end{cases}.
\end{align*}
Since both $\omega_{d_1}$ and $c_{d_2}$ exhibit asymmetry in the way of $\omega_{d_1}(j) = -\omega_{d_1}(n+1-j)$ and $c_{d_2}(j) = -c_{d_2}(n+1-j),$ we may instead write 
\begin{equation}
    \Cov\left(\hat{X}_\inv^{(d_1)}, X_\des^{(d_2)}\right) = 2\sum_{j=1}^{d_2 \wedge (n - d_2)} \omega_{d_1}(j)c_{d_2}(j) = \frac{1}{6} \sum_{j=1}^{d_2 \wedge (n - d_2)} \omega_{d_1}(j)\,. \label{zerocoefficients}
\end{equation}
We now distinguish the cases $d_2 \leq n/2$ and $d_2 > n/2,$ even though this later proves to be redundant. We analyze the range of zero and nonzero coefficients $\omega_{d_1}(j)c_{d_2}(j)$ appearing in the sum in \eqref{zerocoefficients}.

\textbf{Case 1:} $d_2 \leq n/2$ and $d_1 \leq d_2$. Then, $c_{d_2}(j) = 0$ implies $\omega_{d_1}(j) = 0$ (figuratively, the range of zero coefficients solely depends on $d_1$), and 
\begin{align*}
    \Cov\left(\hat{X}_\inv^{(d_1)}, X_\des^{(d_2)}\right) &= \frac{1}{6} \sum_{j=1}^{d_1} (d_1 + 1 - j) = \frac{1}{6}\sum_{j=1}^{d_1} j = \frac{d_1(d_1+1)}{12}\,.
\end{align*}

\textbf{Case 2:} $d_2 \leq n/2$ and $d_1 > d_2$. Then, the range of zero coefficients depends on $d_2$. In light of Figures~\ref{fig4.6} and~\ref{fig4.7}, we also need to distinguish whether $d_2 \leq n-d_1 \Longrightarrow d_1 + d_2 \leq n,$ or $d_1 + d_2 > n$.

\textbf{Case 2a:} $d_2 \leq n/2$ and $d_1 > d_2,$ but $d_1 + d_2 \leq n$. Then,
\begin{align}
    \Cov\left(\hat{X}_\inv^{(d_1)}, X_\des^{(d_2)}\right) &= \frac{1}{6} \sum_{j=1}^{d_2} (d_1 + 1 - j) = \frac{1}{6} \left(\sum_{j=1}^{d_1} (d_1 + 1 - j) - \sum_{j=d_2+1}^{d_1} (d_1 + 1 - j) \right) \nonumber \\
    &= \frac{1}{12}d_1(d_1 + 1) - \frac{1}{6}\Bigl(1 + 2 + \ldots + (d_1 - d_2)\Bigr) \nonumber \\
    &= \frac{1}{12}d_1(d_1 + 1) - \frac{1}{12}(d_1 - d_2)(d_1 - d_2 + 1) \nonumber \\
    &= \frac{1}{6}d_1d_2 - \frac{d_2^2}{12} + \frac{d}{12} = \frac{1}{6}d_1d_2 - \frac{1}{12}d_2(d_2 - 1)\,. \label{case2a}
\end{align}

\textbf{Case 2b:} $d_2 \leq n/2$ and $d_1 + d_2 > n$. Then,
\[
    \Cov\left(\hat{X}_\inv^{(d_1)}, X_\des^{(d_2)}\right) = \frac{1}{6}\sum_{j=1}^{n-d_1} (d_1 + 1 - j) + \frac{1}{6} \sum_{j=n-d_1+1}^{d_2} (n - 2j + 1)\,. 
\]
On one hand, by substituting $d_2 = n-d_1$ in \eqref{case2a}, we have
\begin{subequations}
\begin{align}
    \frac{1}{6}\sum_{j=1}^{n-d_1} (d_1 + 1 - j) &= \frac{1}{6}d_1(n - d_1) - \frac{1}{12}(n - d_1)(n - d_1 + 1) \nonumber \\
    &= \frac{1}{3}nd_1 - \frac{1}{4}d_1^2 - \frac{1}{12}(n^2 - n + d_1)\,. \label{case2b}
\end{align}
On the other hand, we have
\begin{align}
    \frac{1}{6} \sum_{j=n-d_1+1}^{d_2} (n - 2j + 1) &= \frac{1}{6}(n+1)(d_2 + d_1 - n) - \frac{1}{3} \sum_{j=n-d_1+1}^{d_2} j \nonumber \\
    &= \frac{1}{6}(n+1)(d_2 + d_1 - n) - \frac{1}{3} \left(\sum_{j=1}^{d_2} j - \sum_{j=1}^{n-d_1} j \right) \nonumber \\
    &= \frac{1}{6}(n+1)(d_2 + d_1 - n) - \frac{1}{6}d_2(d_2 + 1) + \frac{1}{6}(n - d_1)(n - d_1 + 1) \nonumber \\
    &= \frac{1}{6}\Bigl(d_1^2 - nd_1 - d_2^2 + nd_2\Bigr)\,, \label{case2c}
\end{align}
\end{subequations}
so overall,
\begin{align*}
    \Cov\left(\hat{X}_\inv^{(d_1)}, X_\des^{(d_2)}\right) &= \eqref{case2b} + \eqref{case2c} \\
    &= \frac{1}{6}n(d_1 + d_2) - \frac{1}{6}d_2^2 -  \frac{1}{12}d_1(d_1 + 1) - \frac{1}{12}n(n - 1)\,.
\end{align*}

\textbf{Case 3:} $d_2 > n/2$ and $d_1 \leq n - d_2 \Longleftrightarrow d_1 + d_2 \leq n$. In that case, the range of zero coefficients depends on $d_1,$ and as in Case 1, we obtain 
\[
    \Cov\left(\hat{X}_\inv^{(d_1)}, X_\des^{(d_2)}\right) = \frac{d_1(d_1+1)}{12}\,.
\]

\textbf{Case 4:} $d_2 > n/2$ and $d_1 > n - d_2$. Then, the range of zero coefficients depends on $d_2,$ but to determine whether the structural break of $\omega_{d_1}(j)$ as seen in Figure~\ref{fig4.7} is relevant, we now need to distinguish the subcases $n - d_2 \leq n - d_1 \Longleftrightarrow d_1 \leq d_2$ (no structural break) and $d_1 > d_2$ (structural break). If $d_1 \leq d_2,$ then by replacing $d_2$ with $n - d_2$ in \eqref{case2a}, we get  
\begin{align*}
    \Cov\left(\hat{X}_\inv^{(d_1)}, X_\des^{(d_2)}\right) &= \frac{1}{6}d_1(n-d_2) - \frac{1}{12}(n - d_2)(n - d_2 + 1) \\
    &= \frac{1}{6}n(d_1 + d_2) - \frac{1}{6}d_1d_2 - \frac{1}{12}d_2(d_2 + 1) - \frac{1}{12}n(n-1)\,.
\end{align*}
If $d_1 > d_2,$ then we again replace $d_2$ with $n - d_2$ and obtain the same result as in Case 2b, because \eqref{case2b} is independent of $d_2$ and \eqref{case2c} is invariant under replacing $d_2$ with $n-d_2$. This completes the proof. \qed

\subsection*{Proof of Theorem~\ref{thm3.10}}

To shorten notation, let $\Cov$ be an abbreviation for $\Cov(\hat{X}_\inv^{(d_1)}, X_\des^{(d_2)})$ and let $\sigma_\inv,$ $\sigma_\des$ be abbreviations for $\Var(X_\inv^{(d_1)})^{1/2}$ and $\Var(X_\des^{(d_2)})^{1/2}$, respectively. Recall the notation $\hat{\rho}_n = \Cov / (\sigma_\inv \sigma_\des)$. Moreover, we briefly write $e_2 := n - d_2 = N_{n,d_2}^=$. 

a): Let $n^{1/2} \ll d_1 \ll n$. If $n - d_2 \gg d_1,$ then we have $d_1 + d_2 \leq n$ and therefore $\Cov = \Theta(d_1^2),$ $\sigma_\inv = \Theta(d_1^{3/2})$. If $d_2 \leq n/2,$ we easily find $\hat{\rho}_n = O\left(d_1^2/(d_1^{3/2}\sqrt{n})\right) = \sqrt{d_1}/\sqrt{n} = o(1)$, both when $d_1 \leq d_2$ and when $d_1 > d_2$. If $d_2 > n/2,$ then we have $\sigma_\des = \Theta(\sqrt{e_2}),$ but by the assumption $d_1 = o(e_2)$, we find $\rho_n = o(1)$ from the same arguments. If $n - d_2 \ll d_1,$ i.e., $e_2 = o(d_1)$, we have $d_1 + d_2 > n$. In that case,
\begin{align}
    \Cov\left(X_\inv^{(d_1)}, X_{\des}^{(d_2)}\right) &= \frac{nd_1 + n^2 - ne_2 - d_1n + d_1e_2}{6} - \frac{(n - e_2)(n - e_2 + 1)}{12} - \frac{n(n-1)}{12} \nonumber \\ 
    &= \frac{d_1e_2}{6} - \frac{e_2^2}{12} + \frac{e_2}{12}\,, \label{cov_e2}
\end{align}
yielding
\[
    \hat{\rho}_n = O\left(\frac{d_1e_2 - e_2^2}{d_1^{3/2}\sqrt{e_2}}\right) = O\left(\frac{d_1\sqrt{e_2} - e_2^{3/2}}{d_1^{3/2}}\right) = O\left(\left(\frac{e_2}{d_1}\right)^{1/2} + \left(\frac{e_2}{d_1}\right)^{3/2}\right) = o(1)\,.
\]

b): If $e_2 = cd_1(1 + o(1))$ for some $c > 1$, then $d_1 + d_2 \leq n$, and we now have 
\[
    \hat{\rho}_n = \frac{d_1(d_1 + 1)/12}{\sqrt{(d_1^3/18)(1 + o(1))}\sqrt{cd_1/4}} = \frac{\sqrt{d_1^3}/6}{\sqrt{c}\sqrt{d_1^3/18}}(1 + o(1)) \longrightarrow \frac{1/6}{\sqrt{c/18}} = \frac{1}{\sqrt{2c}}\,.
\]

c): When taking $e_2 = cd_1(1 + o(1))$ with $c < 1$, we have by \eqref{cov_e2}: $\Cov = cd_1^2/6 - c^2d_1^2/12 + cd_1/12 = cd_1(d_1/6 - cd_1/12 + 1/12)$, hence
\[
    \hat{\rho}_n = \frac{cd_1(d_1/6 - cd_1/12 + 1/12)}{\sqrt{d_1^3/18}(1 + o(1))\sqrt{cd_1/4}} = \frac{\sqrt{c}}{d_1}\sqrt{72}\left(\frac{d_1}{6}-\frac{cd_1}{12}+\frac{1}{12}\right) = \sqrt{2c}\left(1 - \frac{c}{2}\right) + o(1)\,.
\]

d): Under the given assumptions $d_1 \asymp c_1n,$ $d_2 \asymp c_2n$ with $c_1 < c_2 < 1/2,$ we have $\Cov = d_1(d_1+1)/12 = c_1^2n^2/12 + O(n)$ and
\[
    \sigma_\inv \sigma_\des = \sqrt{\frac{(c_1n)^3}{18} + O(n^2)}\sqrt{\frac{n + c_2n}{12}} = c_1^{3/2}\frac{n^{3/2}}{\sqrt{18}}(1 + o(1)) \sqrt{(1 + c_2n)/12}\,,
\]
hence
\[
    \hat{\rho}_n \asymp \frac{c_1^2n^2/12}{c_1^{3/2}(n^{3/2}/\sqrt{18})\,\sqrt{n}\sqrt{1+c_2}/\sqrt{12}} = \frac{c_1^2}{c_1^{3/2}\sqrt{1+c_2}} \cdot \frac{1/12}{1/\sqrt{18} \cdot 1/\sqrt{12}} = \frac{\sqrt{6}}{2}\frac{\sqrt{c_1}}{\sqrt{1+c_2}}\,.
\]

e): If $c_1 < 1/2,$ but $1/2 < c_2 \leq 1 - c_1,$ then we have $\sigma_\des = \sqrt{(1 - c_1)n}/2$, while $\Cov$ and $\sigma_\inv$ are the same as in e), giving
\[
    \hat{\rho}_n \asymp \frac{c_1^2n^2/12}{c_1^{3/2}(n^{3/2}/\sqrt{18})\sqrt{(1-c_2)n}/2} = \frac{1}{\sqrt{2}}\frac{\sqrt{c_1}}{\sqrt{1-c_2}}\,.
\]
If $c_2 > 1 - c_1,$ then it follows from \eqref{cov_e2} that $\Cov = (c_1(1 - c_2)/6 - (1 - c_2)^2/12)n^2 + O(n),$ giving 
\[
    \hat{\rho}_n \asymp \frac{c_1(1-c_2)/6 - (1 - c_2)^2/12}{(c_1^{3/2}/\sqrt{18})\,\sqrt{1-c_2}/2} = \frac{\sqrt{1 - c_2}}{\sqrt{2}}\frac{2c_1 + c_2 - 1}{c_1^{3/2}} = \frac{f(c_1,c_2)}{\sqrt{2}c_1^{3/2}}\,.
\]

f): Since now $c_1 > 1/2,$ i.e., $d_1 > n/2,$ we have by Theorem~\ref{thm2.3} that 
\[
    \sigma_\inv = \sqrt{-\frac{1}{6}c_1^3 + \frac{1}{3}c_1^2 - \frac{1}{6}c_1 + \frac{1}{36}}n^{3/2}(1 + o(1)) = \sqrt{\frac{p(c_1)}{3}}n^{3/2}(1 + o(1))\,. 
\]
Then, with the aid of \eqref{cov_e2}, we have
\[
    \hat{\rho}_n \asymp \frac{c_1(1-c_2)/6 - (1 - c_2)^2/12}{\sqrt{p(c_1)/3}\,\sqrt{1-c_2}/2} = \frac{f(c_1, c_2)/12}{\sqrt{p(c_1)}/(2\sqrt{3})} = \frac{\sqrt{3}}{6}\frac{f(c_1, c_2)}{\sqrt{p(c_1)}}\,.
\]

g): If $c_1 > c_2,$ but $c_1 < 1/2$, then $\Cov \asymp (c_1c_2/6 - c_2^2/12)n^2 + O(n),$ giving 
\[
    \hat{\rho}_n \asymp \frac{c_1c_2/6 - c_2^2/12}{c_1^{3/2}/\sqrt{18}\, \sqrt{1+c_2}/\sqrt{12}} = \frac{\sqrt{6}}{2}\frac{2c_1c_2 - c_2^2}{c_1^{3/2}\sqrt{1 + c_2}}\,.
\]

h): If $c_1 > 1/2,$ but $c_2 \leq 1 - c_1$ (forcefully implying $c_2 < n/2$), then compared with g), $\sigma_\inv$ is given by $\sqrt{p(c_1)/3}n^{3/2}(1 + o(1)),$ so that
\[
    \hat{\rho}_n \asymp \frac{c_1c_2/6 - c_2^2/12}{\sqrt{p(c_1)/3}\sqrt{1+c_2}/\sqrt{12}} = \frac{2c_1c_2 - c_2^2}{\sqrt{p(c_1)}\sqrt{1+c_2}\sqrt{12}/\sqrt{3}} = \frac{2c_1c_2 - c_2^2}{2\sqrt{p(c_1)}\sqrt{1+c_2}}\,.
\]

i): If $c_1 > 1/2$ and $c_2 > 1 - c_1$ but $c_2 < 1/2,$ we have $d_1 + d_2 > n,$ hence $\Cov = n(d_1 + d_2)/6 - d_2^2/6 - d_1(d_1 + 1)/12 - n(n-1)/12$ according to Lemma~\ref{lemma3.10}, while still $\sigma_\des = \sqrt{(n + d_2)/12}$. Taking $e_1 = n - d_1 \asymp (1 - c_1)n,$ a similar calculation as in \eqref{cov_e2} gives 
\begin{align*}
    \Cov &= \frac{n((n - e_1) + d_2) - d_2^2}{6} - \frac{(n - e_1)(n - e_1 + 1)}{12} - \frac{n(n-1)}{12} \\
    &= \frac{d_2n}{6} - \frac{d_2^2}{6} - \frac{e_1^2}{12} + O(n) \asymp \frac{c_2(1 - c_2)}{6}n^2 - \frac{(1 - c_1)^2}{12}n^2\,.
\end{align*}
This gives 
\[
    \hat{\rho}_n \asymp \frac{2c_2(1 - c_2) - (1 - c_1)^2}{2\sqrt{p(c_1)}\sqrt{1+c_2}}\,.
\]
Last, if $c_2 > 1/2 > 1 - c_1,$ we also have $\sigma_\des \asymp \sqrt{(1 - c_2)n}/2,$ giving
\[
    \hat{\rho}_n \asymp \frac{(2c_2(1 - c_2) - (1 - c_1)^2)/12}{\sqrt{p(c_1)/3}\, \sqrt{1+c_2}/2} = \frac{\sqrt{3}}{6} \frac{2c_2(1 - c_2) - (1 - c_1)^2}{\sqrt{p(c_1)}\sqrt{1-c_2}}\,. 
\]
This completes the classification and computation of the limiting correlation. \qed

\subsection*{Proof of Lemma~\ref{lemma3.12}} In the proof of Lemma~\ref{lemma3.7}, we established that 
\[
    \hat{X}_\inv^{B,(d_1)} = \sum_{k=1}^n \omega_{d_1}(k)Z_k + \sum_{k=1}^{n \wedge d}\textbf{1}\{Z_k < 0\} + \text{const}\,,
\]
where the coefficients $\omega_{d_1}(k)$ are given in \eqref{omega_a}, \eqref{omega_b}, \eqref{omega_c}, \eqref{omega_d}, depending on whether $d_1 \leq n/2$, or $n/2 < d_1 \leq 2n/3$, or $2n/3 < d_1 < n,$ or $d_1 \geq n$. Similar to the proof of Lemma~\ref{lemma3.10}, we write 
\[
    \Cov\left(\hat{X}_\inv^{B,(d_1)}, X_\des^{B,(d_2)}\right) = \sum_{k=1}^n \omega_{d_1}(k)c_{d_2}(k) + \sum_{k=1}^{n \wedge d} \Cov(\textbf{1}\{Z_k < 0\}, X_\des^{(d_2)})\,,
\]
where $c_{d_2}(k) := \Cov(Z_k, X_\des^{(d_2)})$. 
However, the second sum is clearly $O(d_1)$ and therefore always asymptotically negligible. For any $k$, we have under the assumption that indices are not out of bounds: 
\begin{align*}
    \Cov(Z_k, \textbf{1}\{Z_k > Z_{k+d_2}\}) &= 1/6\,, & \Cov(Z_k, \textbf{1}\{Z_{k-d_2} > Z_k\}) &= -1/6\,, \\
    \Cov(Z_k, \textbf{1}\{-Z_k > Z_{d_2-k}\}) &= -1/6\,, & \Cov(Z_k, \textbf{1}\{-Z_{d_2-k} > Z_k\}) &= -1/6\,. 
\end{align*}
In light of the decomposition of $X_\des^{B,(d_2)}$, we also decompose $c_{d_2}(k) = c_{1}(k) + c_{2}(k) + c_{3}(k)$, 
where
\begin{align*}
    c_{1}(k) &= \Cov\left(Z_k, \sum_{i=1}^d \textbf{1}\{Z_k > Z_{i+d_2}\}\right), \\ 
    c_{2}(k) &= \Cov\left(Z_k, \sum_{i=1}^{\lceil d/2\rceil - 1} \textbf{1}\{-Z_i > Z_{d_2-i}\}\right), \\
    c_{3}(k) &= \Cov\left(Z_k,  \textbf{1}\{Z_{d_2} < 0\}\right). 
\end{align*}
Then, following the lines of the proof of Lemma~\ref{lemma3.10}, we have  
\begin{align*}
    d_2 \leq n/2 &\Longrightarrow c_{1}(k) = \begin{cases} +1/6\,, & k \leq d_2 \\ -1/6\,, & k > n-d_2 \\ 0, & \text{otherwise}  \end{cases}, 
    \\
    d_2 > n/2 &\Longrightarrow c_{1}(k) = \begin{cases} +1/6\,, & k \leq n - d_2 \\ -1/6\,, & k > d_2  \\ 0 & \text{otherwise}  \end{cases}. 
\end{align*} 
In particular, if $d_2 \geq n,$ then $c_{1} \equiv 0$ globally. Next, we have $c_{2}(k) = 0$ for all $k \geq d_2$ and also for $k = d_2/2$ (if $d_2$ is even), since there is no root $(d/2, d/2)$. Further, if $d_2 \geq n$, then $c_2(k) = 0$ also applies for any $k < d_2 - n,$ since the summation of indicators $\textbf{1}\{-Z_i > Z_{d_2-i}\}$ then only spans over $k \in \{d_2 - n, \ldots, \lceil d/2\rceil - 1\}$. In any other case, we have $c_{2}(k) = -1/6$. Finally, if $k = d_2,$ we have $c_{3}(k) = -1/4$, otherwise we have $c_{3}(k) = 0$. Therefore, the influence of $c_{3}(k)$ is negligible, so we rewrite $c_{d_2}(k) = c_1(k) + c_2(k)$. Putting the previous observations together, we have:
\begin{align*}
    d_2 \leq n/2 &\Longrightarrow c_{d_2}(k) = \begin{cases} 0, & k < d_2/2 \\ 1/6, & k=d_2/2~\text{or}~k=d_2 \\ 0, & d_2/2 < k < d_2~\text{or}~d_2 < k \leq n - d_2\\ -1/6, & k > n - d_2 \end{cases}, 
    \\
    n/2 < d_2 \leq 2n/3 &\Longrightarrow c_{d_2}(k) = \begin{cases} 0, & k < d_2/2 \\ 1/6, & k=d_2/2 \\ 0, & d_2/2 < k \leq n-d_2 \\ -1/6, & k > n - d_2 \end{cases}, 
    \\
    2n/3 < d_2 < n &\Longrightarrow c_{d_2}(k) = \begin{cases} 0, & k \leq n - d_2 \\ -1/6, & n - d_2 < k < d_2/2 \\ 0, & k = d_2/2 \\ -1/6, & k > d_2/2 \end{cases}. 
\end{align*} 
For all $d_2 \geq n,$ we simply have 
\[
    c_{d_2}(k) = c_2(k) = \begin{cases} 0, & k < d_2 - n \\ -1/6, & k \geq d_2 - n, k \neq d_2/2 \\ 0, & k = d_2/2 \end{cases}.
\]
In conclusion, the identity
\[
    c_{d_2}(k) = \begin{cases} 0, & k \leq |n - d_2| \\ -1/6, & k > |n - d_2| \end{cases}
\]
applies for all $k$ with at most two exceptions ($k = d_2/2$ and $k = d_2).$ These exceptions are asymptotically negligible throughout. We likewise find that if $d_1 < n,$ then
\[
    \omega_{d_1}(k) = \begin{cases}
        O(1), & k \leq n - d_1 \\ (n - d_1 - k)/2 + O(1), & k > n - d_1
    \end{cases}.
\]
If $d_1 > n,$ then all $k = 1, \ldots, n$ have to be taken into account, since by \eqref{omega_d}, 
\[
    \omega_{d_1}(k) = \begin{cases}
        1 - k, & k \leq d_1 - n \\ (n - d_1 - k)/2 + O(1), & k > d_1 - n
    \end{cases}.
\]
The $O(1)$ remainders for $k > |n-d_1|$ are clearly dominated by the $(n-d_1-k)/2$ entries and only add $O(d_1)$ overall. Regarding the $O(1)$ entries for $d_1 < n$ and $k \leq n-d_1$, we find in \eqref{omega_a} that $\omega_{d_1}(k) = 0$ holds if $d_1 < k \leq n - d_1$. Especially, for sublinear regimes $d_1 = o(n),$ all but at most $2d_1$ coefficients vanish altogether, so the $O(1)$ entries for $d_1 < n$ also contribute no more than $O(d_1)$. We can therefore focus on the above simplifications of $\omega_{d_1}(k)$ and $c_{d_2}(k)$, and distinguish the following cases: 

\textbf{Case 1:} $d_1 < n$ and $d_2 < n$. Then, the relevant indices are those satisfying $k > \max\{n - d_1, n - d_2\}$ and we need to distinguish whether $d_1 \geq d_2$ or $d_1 < d_2$. In the former case, we have 
\begin{align*}
    \Cov\left(\hat{X}_\inv^{B,(d_1)}, X_\des^{B,(d_2)}\right) &= -\frac{1}{12}\sum_{k=n-d_2+1}^n (n - d_1 - k) \\
    &= -\frac{1}{12} \sum_{k=n-d_2+1} (n - d_1 - d_2 - k + d_2) \\
    &= -\frac{1}{12} \sum_{k=n-d_2+1} (n - d_2 - k) - \frac{1}{12} d_2(d_2 - d_1) \\
    &= \frac{1}{24}d_2(d_2+1) + \frac{1}{12}d_2(d_1 - d_2)\,. 
\end{align*}
In the latter case of $d_1 < d_2,$ we simply have 
\[
    \Cov\left(\hat{X}_\inv^{B,(d_1)}, X_\des^{B,(d_2)}\right) = -\frac{1}{12}\sum_{k=n-d_1+1}^n (n - d_1 - k) = \frac{d_1(d_1+1)}{24}\,.
\]

\textbf{Case 2:} $d_1 > n,$ but $d_2 < n$. Since there is a structural break for $\omega_{d_1}(k)$ at $k = d_1 - n,$ we now need to compare $d_1 - n$ with $n - d_2,$ i.e., distinguish whether $d_1 + d_2 > 2n$ or $d_1 + d_2 \leq 2n$. In the former case, we have
\begin{align*}
    &\Cov\left(\hat{X}_\inv^{B,(d_1)}, X_\des^{B,(d_2)}\right) \\
    =~& -\frac{1}{6}\sum_{k=n-d_2+1}^{d_1-n} (1-k) - \frac{1}{6}\sum_{k=d_1-n+1}^n (n-d_1-k)/2 \\
    =~& -\frac{1}{6}(d_1 + d_2 - 2n) + \frac{1}{12}(d_1 - n)(d_1 - n + 1) - \frac{1}{12}(n - d_2)(n - d_2 + 1) \\
    & \quad - \frac{1}{12}\sum_{k=d_1-n+1}^n (n - d_1) + \frac{1}{12}\sum_{k=d_1-n+1}^n k \\
    =~& -\frac{1}{6}(d_1 + d_2 - 2n) + \frac{1}{12}(d_1 - n)(d_1 - n + 1) - \frac{1}{12}(n - d_2)(n - d_2 + 1) \\
    & \quad -\frac{1}{12}(2n - d_1)(n-d_1) + \frac{1}{24}n(n+1) - \frac{1}{24}(d_1 - n)(d_1 - n + 1) \\
    =~& -\frac{1}{24}d_1(d_1 - 4n + 3) - \frac{1}{12}d_2(d_2 - 2n + 1) - \frac{n^2}{6} + \frac{n}{4}\,.
\end{align*}
If $d_1 + d_2 \leq 2n,$ then the calculation simplifies to
\begin{align*}
    \Cov\left(\hat{X}_\inv^{B,(d_1)}, X_\des^{B,(d_2)}\right) &= -\frac{1}{6} \sum_{k=n-d_2+1}^n (n - d_1 - k)/2 \\
    &= \frac{1}{24}d_2(d_2+1) + \frac{1}{12}d_2(d_1 - d_2)\,.
\end{align*}

\pagebreak 
\textbf{Case 3:} If $d_1 < n$ but $d_2 > n,$ then the interfering structural breaks for $\omega_{d_1}(k)$ and $c_{d_2}(k)$ are given by $n - d_1$ and $d_2 - n,$ respectively, which means we again compare the cases $d_1 + d_2 \leq 2n$ and $d_1 + d_2 > 2n$. In the former case, we simply have
\[
    \Cov\left(\hat{X}_\inv^{B,(d_1)}, X_\des^{B,(d_2)}\right) = -\frac{1}{12}\sum_{k=n-d_1+1}^n (n - d_1 - k) = \frac{d_1(d_1+1)}{24}\,.
\]
In the latter case, we calculate 
\begin{align*}
    \Cov\left(\hat{X}_\inv^{B,(d_1)}, X_\des^{B,(d_2)}\right) &= -\frac{1}{12} \sum_{k=d_2-n}^n (n - d_1 - k) \\
    &= -\frac{1}{12} \left(\sum_{k=d_2-n}^n (n - k) - d_1(2n - d_2 + 1)\right) \\
    &= -\frac{1}{12}\left(\frac{(2n - d_2)(2n - d_2 + 1)}{2} - d_1(2n - d_2 + 1)\right) \\
    &= \frac{1}{12}(2n - d_2 + 1)(d_1 - n + d_2/2)\,.
\end{align*}

\textbf{Case 4:} If both $d_1 > n$ and $d_2 > n,$ then in case of $d_1 \leq d_2,$ the structural breaks behave the same as in the latter subcase of Case 3, i.e.,
\[
    \Cov\left(\hat{X}_\inv^{B,(d_1)}, X_\des^{B,(d_2)}\right) = -\frac{1}{12} \sum_{k=d_2-n}^n (n - d_1 - k) = \frac{1}{12}(2n - d_2 + 1)(d_1 - n + d_2/2)\,.
\]
Otherwise, i.e., if $d_1 > d_2,$ we proceed similarly as in Case 2, obtaining
\begin{align*}
    &\Cov\left(\hat{X}_\inv^{B,(d_1)}, X_\des^{B,(d_2)}\right) \\
    =~& -\frac{1}{6}\sum_{k=d_2-n+1}^{d_1-n} (1-k) - \frac{1}{6}\sum_{k=d_1-n+1}^n (n-d_1-k)/2 \\
    =~& -\frac{1}{6}(d_1 - n - (d_2 - n + 1) + 1) + \frac{1}{6}\sum_{k=d_2-n+1}^{d_1-n} k - \frac{1}{12}(2n - d_1)\left(n - \frac{3}{2}d_1 - \frac{1}{2}\right) \\
    =~& -\frac{1}{6}(d_1 - d_2) + \frac{1}{12}(d_1 - n)(d_1 - n + 1) - \frac{1}{12}(d_2 - n)(d_2 - n + 1) \\
    & \quad - \frac{1}{12}(2n - d_1)\left(n - \frac{3}{2}d_1 - \frac{1}{2}\right) \\
    =~& -\frac{d_1}{24}(d_1 - 4n + 3) - \frac{1}{12}d_2(d_2 - 2n - 1) - \frac{n^2}{6} + \frac{n}{12}\,.
\end{align*}
This completes the calculation of $\Cov\left(\hat{X}_\inv^{B,(d_1)}, X_\des^{B,(d_2)}\right)$. \qed

\section{Calculations of original covariances and correlations} \label{appendixA}

\noindent \textbf{Proof of Theorem~\ref{thm3.13}}. On symmetric groups, the covariance of $X_\inv^{(d_1)}$ and $X_\des^{(d_2)}$ is 
\begin{align}
    \Cov(X_\inv^{(d_1)}, X_\des^{(d_2)}) &= \sum_{(i,j)\in \NN_{n,d_1}}\sum_{k=1}^{n-d_2} \Cov(\textbf{1}\{Z_i > Z_j\}, \textbf{1}\{Z_k > Z_{k+d_2}\}) \nonumber \\
    &= \sum_{(i,j) \in \NN_{n,d_1}} S_{ij}\,, \label{sij}
\end{align}
taking
\[
    S_{ij} := \sum_{k=1}^{n-d_2} \Cov(\textbf{1}\{Z_i > Z_j\}, \textbf{1}\{Z_k > Z_{k+d_2}\})\,.
\]
Consider an arbitrary pair $(i,j) \in \NN_{n,d_1}$. For any $k \notin \{i-d_2, i, j-d_2, j\},$ the events $\{Z_i > Z_j\}$ and $\{Z_k > Z_{k+d_2}\}$ are independent. Hence, the relevant \textit{candidate indices} are $k=i-d_2, i, j-d_2, j,$ and we can write
\begin{align*}
    S_{ij} &= \sum_{k \in \{i-d_2, i, j-d_2, j\} \cap \{1,\ldots,n-d_2\}} \Cov(\textbf{1}\{Z_i > Z_j\}, \textbf{1}\{Z_k > Z_{k+d_2}\}) \\
    &= \sum_{k \in \{i-d_2, i, j-d_2, j\}\cap\{1,\ldots,n-d_2\}} \Bigl[\,\PP(Z_i > Z_j, Z_k > Z_{k+d_2}) - 1/4\,\Bigr]\,.
\end{align*}
To compute each $S_{ij},$ we first distinguish the following cases:
\begin{itemize}
    \item If $i-j < d_2,$ then the order of $k$-candidates is given by $i-d_2 < j-d_2 < i < j$. We call this an \textit{intersection} (since the intersection of $\{i-d_2, \ldots, i\}$ and $\{j-d_2,\ldots,j\}$ is nonempty). 
    \item If $i-j = d_2,$ i.e., $i=j-d_2,$ then we have a \textit{collision} (as the candidates $k=i$ and $k=j-d_2$ collide). This is only possible if $d_1 \geq d_2$.
    \item If $i-j > d_2,$ then the order of $k$-candidates is given by $i-d_2 < i < j-d_2 < j$. We call this a \textit{separation}. This is only possible if $d_1 > d_2$.
\end{itemize}
In conclusion, we decompose the covariance into three parts, namely, 
\[
    \sum_{(i,j) \in \NN_{n,d_1} \cap \NN_{n,d_2}^<} S_{ij} + \sum_{(i,j) \in \NN_{n,d_1} \cap \NN_{n,d_2}^=} S_{ij} + \sum_{(i,j) \in \NN_{n,d_1} \setminus \NN_{n,d_2}} S_{ij}\,.
\]
At first, assume $i \neq j - d_2$. For $k=i-d_2$ and $k=j,$ we get $\PP(Z_i > Z_j, Z_k > Z_{k+d_2}) = 1/6,$ while for $k=i$ and $k=j-d_2,$ we get $\PP(Z_i > Z_j, Z_k > Z_{k+d_2}) = 1/3.$ The resulting contributions to the covariance are 
\[
    \Cov(\textbf{1}\{Z_i > Z_j\}, \textbf{1}\{Z_k > Z_{k+d_2}\}) = \begin{cases} 1/12, & k=j-d_2~\text{or}~k=i \\ -1/12, & k=i-d_2~\text{or}~k=j \end{cases}.
\]

If $\{i-d_2, i, j-d_2, j\} \subseteq \{1,\ldots, n-d_2\}$, then these four contributions cancel to zero. Therefore, we need to analyze the constellations $(i,j)$ for which at least one of the candidates $k=i-d_2,i,j-d_2,j$ is outside the permitted index set $\{1, \ldots, n-d_2\}$. Moreover, in the collision case $i=j-d_2,$ we have only three candidates $k \in \{i-d_2, i, i+d_2\},$ while
\[
    \Cov(\textbf{1}\{Z_i > Z_{i+d_2}\}, \textbf{1}\{Z_k > Z_{k+d_2}\}) = \begin{cases} 1/4, & k=i \\ -1/12, & k=i-d_2~\text{or}~k=i+d_2\end{cases}.
\]
The non-colliding candidates $k=i-d_2$ and $k=j$ may also possibly be outside the permitted index set. In what follows, we briefly say that the candidates outside $\{1, \ldots, n-d_2\}$ \textit{vanish}, while those staying inside $\{1, \ldots, n-d_2\}$ \textit{survive}. In general, we can note the following: 
\begin{itemize}
    \item $k=i-d_2$ or $k=j-d_2$ vanishes iff $i \leq d_2$ and $j \leq d_2,$ respectively.
    \item $k=i$ or $k=j$ vanishes iff $i > n-d_2$ or $j > n-d_2,$ respectively.
\end{itemize} 

By the number of surviving candidates, we categorize the contributions $S_{ij}$ as follows:

\noindent \textbf{Intersection} $(j-i<d_2)$:
\begin{itemize}
    \item \textit{three candidates:} only $i-d_2$ vanishes, or only $j$ vanishes, giving $S_{ij} = 1/12$. These cases are symmetric by interchanging $i'=n+1-j$ and $j'=n+1-i$. The case of only $i$ vanishing (or, symmetrically, only $j-d_2$ vanishing) is not possible, as the requirements $i > n-d_2$ and $j \leq n - d_2$ are contradictory.
    \item \textit{two candidates:} if both $i-d_2$ and $j-d_2$ vanish, or if both $i$ and $j$ vanish, we still have $S_{ij} = 1/12 - 1/12 = 0$. However, if both $i-d_2$ and $j$ vanish (while $j-d_2$ and $i$ survive), we have $S_{ij} = 1/12 + 1/12 = 1/6$. 
    \item \textit{one candidate:} It is not possible that one of $k=i-d_2$ or $k=j$ is the only surviving candidate, since this again yields contradictions. The case of $k=j-d_2$ or $k=i$ being the only surviving candidate is possible if $d_2$ is sufficiently large, and such constellations contribute $S_{ij} = 1/12$.
    \item Cases in which all four candidates vanish need not be analyzed.
\end{itemize}

\noindent \textbf{Collision} $(j-i=d_2)$:
\begin{itemize}
    \item \textit{all three candidates:} $k=i$ always survives. If all candidates $k=i-d_2, i, i+d_2$ survive, then we have $S_{ij} = 1/4 - 1/12 - 1/12 = 1/12$.
    \item \textit{two candidates:} If only $k=i-d_2$ or only $k=j=i+d_2$ vanishes, then we have $S_{ij} = 1/4 - 1/12 = 1/6$.
    \item \textit{one candidate:} If both $k=i-d_2$ and $k=j=i+d_2$ vanish, we have  $S_{ij} = 1/4$. 
\end{itemize}

\noindent \textbf{Separation} $(j-i>d_2)$: The cases of three or two candidates are categorized in the same way as for intersections. However, cases of only one candidate surviving are not possible, since they require $i > n-d_2$ or $j < d_2,$ in contradiction to the assumption $j-i>d_2$. \par
\medskip
Taking symmetric constellations into account, these considerations are formalized as follows:
\begin{align*}
    C_1 =~& \frac{1}{6}\#\{(i,j) \in \NN_{n,d_1} \cap \NN_{n,d_2}^< : i \leq d_2, i\leq n-d_2, j>d_2, j \leq n-d_2\} \\
    &+ \frac{1}{6}\#\{(i,j) \in \NN_{n,d_1} \cap \NN_{n,d_2}^< : i \leq d_2, i \leq n - d_2, j > d_2,  j > n-d_2\} \\
    &+ \frac{1}{6}\#\{(i,j) \in \NN_{n,d_1} \cap \NN_{n,d_2}^< : i \leq d_2, i > n -d_2, j > d_2, j > n-d_2\} \\
    =:~&C_{11} + C_{12} + C_{13}\,, \\
    C_2 =~& \frac{1}{4}\#\{i = 1,\ldots,n-d_2: i \leq d_2, i > n - 2d_2\} \\
    &+ \frac{1}{3}\#\{i = 1,\ldots,n-d_2: i\leq d_2, i \leq n-2d_2\} \\
    &+ \frac{1}{12}\#\{i = 1,\ldots,n-d_2: i > d_2, i \leq n - 2d_2\} \\
    =:~& C_{21} + C_{22} + C_{23}\,, \\
    C_3 =~& \frac{1}{6}\#\{(i,j) \in \NN_{n,d_1} \setminus \NN_{n,d_2} : i \leq d_2, i \leq n-d_2, j>d_2, j \leq n-d_2\} \\
    &+ \frac{1}{6}\#\{(i,j) \in \NN_{n,d_1} \setminus \NN_{n,d_2} : i \leq d_2, i \leq n - d_2, j > d_2, j > n-d_2\} \\
    =:~& C_{31} + C_{32}\,,
\end{align*}
such that
\[
    \Cov(X_\inv^{(d_1)}, X_\des^{(d_2)}) = C_1 + C_2 \cdot \textbf{1}\{d_1 \geq d_2\} + C_3\,.
\]
We now need to distinguish cases as follows: 
\begin{itemize}
\item If $d_2 < n/2,$ then $C_{13} = 0$ and 
\begin{align*}
    C_{11} &= \frac{1}{6}\#\{(i,j) \in \NN_{n,d_1} \cap \NN_{n,d_2}^< : i \leq d_2, d_2 < j \leq n - d_2\}\,, \\
    C_{12} &= \frac{1}{6}\#\{(i,j) \in \NN_{n,d_1} \cap \NN_{n,d_2}^< : i \leq d_2, j > n - d_2\}\,,
\end{align*}
yielding $C_{11} + C_{12} = \displaystyle{\frac{1}{6}}\#\{(i,j) \in \NN_{n,d_1} \cap \NN_{n,d_2}^<: i \leq d_2, j > d_2\}.$
\item If $d_2 = n/2,$ then $C_{11} = C_{13} = 0$ and
\[
    C_{12} = \frac{1}{6}\#\{(i,j) \in \NN_{n,d_1} \cap \NN_{n,d_2}^< : i \leq n/2, j > n/2\}\,.
\]
\item If $d_2 > n/2,$ then $C_{11} = 0$ and 
\begin{align*}
    C_{12} &= \frac{1}{6}\#\{(i,j) \in \NN_{n,d_1} \cap \NN_{n,d_2}^< : i \leq n - d_2, j > d_2\}\,, \\
    C_{13} &= \frac{1}{6}\#\{(i,j) \in \NN_{n,d_1} \cap \NN_{n,d_2}^< : n - d_2 < i \leq d_2, j > d_2\}\,,
\end{align*}
yielding $C_{12} + C_{13} = \displaystyle{\frac{1}{6}}\#\{(i,j) \in \NN_{n,d_1} \cap \NN_{n,d_2}^<: i \leq d_2, j > d_2\}.$
\end{itemize}
Therefore, we have in any case: 
\[
    C_1 = \frac{1}{6}\#\{(i,j) \in \NN_{n,d_1} \cap \NN_{n,d_2}^< : i \leq d_2, j > d_2\}\,. 
\]
Without restriction, assume $d_1 < d_2$ (otherwise, we set $d_1 = d_2 - 1$ for now). First, consider the case $d_1 + d_2 \leq n,$ so that for every index $i \in \{1, \ldots, d_2\}$, we have $i + d_1 \leq n$ (meaning that all pairs in $\NN_{n,d_1}$ involving $i$ as first index can be taken into account). Then we have
\begin{align*}
    C_1 &= \frac{1}{6}\#\{(i,j) \in \NN_{n,d_1} : i \leq d_2, j > d_2\} \\
    &= \frac{1}{6}\Bigl(\#\{(d_2, d_2+1), \ldots, (d_2, d_2+d_1)\} \\
    & \qquad \quad + \#\{(d_2 - 1, d_2 + 1), \ldots, (d_2 - 1, d_2 + d_1 - 1)\} \\
    & \qquad \quad + \ldots + \#\{(d_2 - d_1 + 1, d_2 + 1)\}\Bigr) = \frac{d_1(d_1 + 1)}{12}.
\end{align*}
In the case of $d_1 + d_2 > n,$ we have to subtract the number of pairs $(i,j)$ for which $i \leq d_2$ and $j - i \leq d_1,$ but $j > n$. Therefore, we now obtain
\begin{align*}
    C_1 &= \frac{1}{6}\left(\frac{d_1(d_1 + 1)}{2} - \frac{(d_1 + d_2 - n)(d_1 + d_2 - n + 1)}{2}\right) \\
    &= \frac{1}{6}\left(n(d_1 + d_2) - d_1d_2\right) - \frac{1}{12}\left(d_2(d_2 + 1) + n(n-1)\right).
\end{align*}
In the case of $d_1 < d_2$ and $d_1 + d_2 \leq n,$ we also have the final result
\[
    \Cov(X_\inv^{(d_1)}, X_\des^{(d_2)}) = \frac{d_1(d_1+1)}{12}\,.
\]
In the case of $d_1 \geq d_2,$ plugging $d_1 = d_2 - 1$ yields 
\[
    C_1 = \begin{cases} d_2(d_2 - 1)/12, & d_2 \leq (n+1)/2 \\ \frac{1}{3}d_2n - \frac{1}{4}d_2^2 + \frac{1}{12}(d_2 - n(n+1)), & d_2 > (n+1)/2 \end{cases}.
\]

Next, we turn to $C_3$. Recall that $C_{3} = C_{31} + C_{32}$. By previous arguments, we find analogous representations for $C_{31}$ and $C_{32},$ namely:
\begin{align*}
    C_{31} &= \begin{cases} \#\{(i,j) \in \NN_{n,d_1} \setminus \NN_{n,d_2}: i \leq d_2, d_2 < j \leq n-d_2\}/6, & d_2 < n/2 \\ 0, & d_2 \geq n/2 \end{cases}, \\
    C_{32} &= \begin{cases} \#\{(i,j) \in \NN_{n,d_1} \setminus \NN_{n,d_2}: i \leq d_2, j > n-d_2\}/6, & d_2 < n/2 \\ \#\{(i,j) \in \NN_{n,d_1} \setminus \NN_{n,d_2}: i \leq n-d_2, j > d_2\}/6, & d_2 \geq n/2 \end{cases}.
\end{align*}
For any $(i,j) \in \NN_{n,d_1} \setminus \NN_{n,d_2}$, the condition $j > d_2$ is redundant, which gives
\begin{align}
    C_3 &= \begin{cases} \#\{(i,j) \in \NN_{n,d_1} \setminus \NN_{n,d_2}: i \leq d_2\}/6, & d_2 \leq n/2 \\ \#\{(i,j) \in \NN_{n,d_1} \setminus \NN_{n,d_2}: i \leq n-d_2\}/6, & d_2 \geq n/2 \end{cases} \nonumber \\
    &= \frac{1}{6}\#\{(i,j) \in \NN_{n,d_1} \setminus \NN_{n,d_2}: i \leq d_2 \wedge (n - d_2)\}\,. \label{c3}
\end{align} 
By assuming a low regime for $d_1,$ we also have $d_2 < n/2$ and $d_1 + d_2 < n$. Then, we can easily write 
\[
    C_3 = \frac{1}{6}\#\bigcup_{i=1}^{d_2} \{(i, i + d_2 + 1), \ldots, (i, i + d_1 - 1), (i, i + d_1)\} = \frac{d_2(d_1 - d_2)}{6}\,.
\]
Regarding the collisions, we can also use that $d_2 < n/3$, which then gives $C_{21} = 0$, $C_{22} = \#\{i: i \leq d_2\}/3 = d_2/3$ and $C_{23} = \#\{i: i > d_2, i \leq n - 2d_2\}/12 = (n - 3d_2)/12$.  Overall, $C_2 = d_2/3 + (n-3d_2)/12$. 
By taking appropriate case distinctions and summing $C_1,$ $C_2$ and $C_3$, the claim follows. \qed
\par \bigskip \noindent
\textbf{Proof of Theorem~\ref{thm3.14}.} If $d_1 = o\left(n^{1/3}\log(n^{4/3})^{-7/3}\right)$, then by Theorem~\ref{thm2.3}, 
\[
    \Var(X_\inv) = \Theta(d_1n)\,.
\]
If $d_1 < d_2$ and $d_1 + d_2 < n,$ then $\Cov\left(X_\inv^{(d_1)}, X_{\des}^{(d_2)}\right) = d_1(d_1 + 1)/12 = \Theta(d_1^2),$ which is clearly dominated by $\sigma_\inv \sigma_\des = \Theta(\sqrt{d_1 n} \sqrt{n - d_2}),$ regardless the choice of $d_2,$ so $\rho = 0$. If $d_1 + d_2 > n,$ then we can use \eqref{cov_e2} (with $e_2 := n - d_2$) since the covariance at hand is the same as that of $\hat{X}_\inv^{(d_1)}$ and $X_\des^{(d_2)}$. As $\sigma_\des = \Theta(e_2)$ and $e_2 < d_1 = o(n^{1/3}),$ we then have
\[
    \corr\left(X_\inv^{(d_1)}, X_{\des}^{(d_2)}\right) = \Theta\left(\frac{d_1e_2 - e_2^2 + e_2}{\sqrt{nd_1e_2}}\right) \lesssim \frac{\sqrt{d_1e_2} - e_2}{\sqrt{n}} = \frac{o(n^{1/3})}{\sqrt{n}} = o(1)\,.
\]
This proves a). The interesting case arises if $d_1 \geq d_2,$ and $d_1$ eventually remains constant as $n\rightarrow \infty$. By Theorem~\ref{thm3.13}, we have $\Cov\left(X_\inv^{(d_1)}, X_{\des}^{(d_2)}\right) = (n + d_1^2)/12 = n(1 + o(1))/12$ if $d_1 = d_2,$ and also $\Cov\left(X_\inv^{(d_1)}, X_{\des}^{(d_2)}\right) = (n - d_2^2)/12 + d_1d_2/6 = n(1 + o(1))/12$ if $d_1 > d_2$. Then,   
\[
    \corr\left(X_\inv^{(d_1)}, X_{\des}^{(d_2)}\right) = \frac{n(1 + o(1))/12}{\sqrt{nd_1/12 + o(n)}\sqrt{n/12}} = \frac{1}{\sqrt{d_1}} + o(1)\,.
\]
This proves b), and c) follows from the same argument. \qed
\par \bigskip
\noindent \textbf{Proof of Theorem~\ref{thm3.15}.} Again, we carry out this proof only for signed permutation groups, as the arguments work analogously on even-signed permutation groups. Recall that by Theorem~\ref{thm2.8}, the magnitudes $\sigma_\inv = \Theta(\sqrt{nd_1})$ and $\sigma_\des = \Theta(\sqrt{2n-d_2})$ apply on signed permutation groups as well. Therefore, $\rho \neq 0$ can only be achieved if $\Cov\left(X_\inv^{B,(d_1)}, X_\des^{B,(d_2)}\right) = \Theta(n),$ or $\Cov\left(X_\inv^{B,(d_1)}, X_\des^{B,(d_2)}\right) = \Theta(\sqrt{n}e_2),$ in case $e_2 = 2n - d_2$ is sublinear. According to  \eqref{2.7a} and \eqref{2.7b}, both $X_\inv^{B,(d_1)}$ and $X_\des^{B,(d_2)}$ consist of three parts each, so that $\Cov\left(X_\inv^{B,(d_1)}, X_\des^{B,(d_2)}\right)$ is decomposed into nine parts:
\begin{align}
    &\Cov\left(X_\inv^{B,(d_1)}, X_\des^{B,(d_2)}\right) \nonumber \\
    =~& \Cov\left(\sum_{j-i \leq d_1} \textbf{1}\{Z_i > Z_j\}, \sum_{i=1}^{n-d_2} \textbf{1}\{Z_i > Z_{i+d_2}\}\right) \tag{A1} \label{A1} \\
    +~& \Cov\left(\sum_{i+j \leq d_1} \textbf{1}\{-Z_i > Z_j\}, \sum_{i=1}^{n-d_2} \textbf{1}\{Z_i > Z_{i+d_2}\}\right) \tag{A2} \label{A2} \\
    +~& \Cov\left(\sum_{i=1}^{n \wedge d_1} \textbf{1}\{Z_i < 0\}, \sum_{i=1}^{n-d_2} \textbf{1}\{Z_i > Z_{i+d_2}\}\right) \tag{A3} \label{A3} \\
    +~& \Cov\left(\sum_{j-i \leq d_1} \textbf{1}\{Z_i > Z_j\}, \sum_{i=1}^{\lceil d_2/2\rceil - 1} \textbf{1}\{-Z_i > Z_{d_2-i}\}\right) \tag{B1} \label{B1} \\
    +~& \Cov\left(\sum_{i+j \leq d_1} \textbf{1}\{-Z_i > Z_j\}, \sum_{i=1}^{\lceil d_2/2\rceil - 1} \textbf{1}\{-Z_i > Z_{d_2-i}\}\right) \tag{B2} \label{B2} \\
    +~& \Cov\left(\sum_{i=1}^{n \wedge d_1} \textbf{1}\{Z_i < 0\}, \sum_{i=1}^{\lceil d_2/2\rceil - 1} \textbf{1}\{-Z_i > Z_{d_2-i}\}\right) \tag{B3} \label{B3} \\
    +~& \Cov\left(\sum_{j-i \leq d_1} \textbf{1}\{Z_i > Z_j\}, \textbf{1}\{Z_d < 0\}\right) \tag{C1} \label{C1} \\
    +~& \Cov\left(\sum_{i+j \leq d_1} \textbf{1}\{-Z_i > Z_j\}, \textbf{1}\{Z_d < 0\}\right) \tag{C2} \label{C2} \\
    +~& \Cov\left(\sum_{i=1}^{n \wedge d_1} \textbf{1}\{Z_i < 0\}, \textbf{1}\{Z_d < 0\}\right). \tag{C3} \label{C3}
\end{align}
\eqref{A1} is computed by the same line as in Theorem~\ref{thm3.13}, since it is based on a sum of probabilities of the forms $\PP(Z_i > Z_j)$, $\PP(Z_i > Z_j > Z_k)$, and $\PP(Z_i > Z_j, Z_i > Z_k)$,  which remain the same if $Z \sim U(-1,1)$ is employed instead of $Z \sim U(0,1)$. Therefore, \eqref{A1} inherits its behavior from Theorem~\ref{thm3.14}. We will further show that all other parts cannot contribute more than $O(d_1^2)$ nonzero terms, which means they prove insignificant to the limiting correlation in any case. For the parts \eqref{A3}, \eqref{B3}, \eqref{C1}, \eqref{C2}, and \eqref{C3}, this is easy to see: \eqref{C3} is one term at best, and  
\begin{align*}
    \eqref{C1} =~& \sum_{(i, d_2) \in \mathfrak{N}_{n,d_1}} \Cov(\textbf{1}\{Z_i > Z_{d_2}\}, \textbf{1}\{Z_{d_2} < 0\}) \\
    & \qquad + \sum_{(d_2, j) \in \mathfrak{N}_{n,d_1}} \Cov(\textbf{1}\{Z_{d_2} > Z_j\}, \textbf{1}\{Z_{d_2} < 0\})\,, \\
    \eqref{C2} =~& \sum_{(i, d_2) \in \widetilde{\mathfrak{N}}_{n,d_1}} \Cov(\textbf{1}\{-Z_i > Z_{d_2}\}, \textbf{1}\{Z_{d_2} < 0\}) \\
    & \qquad + \sum_{(d_2, j) \in \widetilde{\mathfrak{N}}_{n,d_1}} \Cov(\textbf{1}\{-Z_{d_2} > Z_j\}, \textbf{1}\{Z_{d_2} < 0\}) \\
\end{align*}
have $O(d_1)$ terms each, and \eqref{A3}, \eqref{B3} follow the same structure. 

Regarding \eqref{A2}, we decompose, by analogy with \eqref{sij}:
\begin{align*}
    \eqref{A2} &= \sum_{i+j\leq d_1} \sum_{k=1}^{n-d_2} \Cov(\textbf{1}\{-Z_i > Z_j\}, \textbf{1}\{Z_k > Z_{k+d_2}\}) =: \sum_{(i,j) \in \widetilde{\mathfrak{N}}_{n,d_1}} S_{ij}\,, \\
    S_{ij} &= \sum_{k=1}^{n-d_2} \Cov(\textbf{1}\{-Z_i > Z_j\}, \textbf{1}\{Z_k > Z_{k+d_2}\}) \\
    &= \Cov(\textbf{1}\{-Z_i > Z_j\}, \textbf{1}\{Z_{i-d_2}>Z_i\}) + \Cov(\textbf{1}\{-Z_i > Z_j\}, \textbf{1}\{Z_{i}>Z_{i+d_2}\}) \\
    &+ \Cov(\textbf{1}\{-Z_i > Z_j\}, \textbf{1}\{Z_{j-d_2}>Z_j\}) + \Cov(\textbf{1}\{-Z_i > Z_j\}, \textbf{1}\{Z_{j}>Z_{j+d_2}\})\,,
\end{align*}
and we already see that $S_{ij}$ encompasses four terms at most (which additionally exhibit cancellation), and $\widetilde{\mathfrak{N}}_{n,d_1}$ is of size $O(d_1^2)$. The same consideration applies for \eqref{B2}. Finally, consider \eqref{B1} in the same fashion:
\begin{align*}
    \eqref{B1} &= \sum_{j-i \leq d_1} S_{ij}, & 
    S_{ij} &:= \sum_{k=1}^{\lceil d_2/2\rceil - 1} \Cov\left( \textbf{1}\{Z_i > Z_j\},  \textbf{1}\{-Z_k > Z_{d_2-k}\}\right).
\end{align*}
Obviously, $S_{ij} = 0$ if $i, j > d_2$. We now consider the case of $i, j \leq d_2$. Then, barring the exceptional case $i + j = d_2,$ we have
\begin{align*}
    S_{ij} &= \Cov\left( \textbf{1}\{Z_i > Z_j\}, \textbf{1}\{-Z_i > Z_{d_2-i}\}\right) + \Cov\left( \textbf{1}\{Z_i > Z_j\},  \textbf{1}\{-Z_j > Z_{d_2-j}\}\right) \\
    &= \Cov\left( \textbf{1}\{Z_1 > Z_2\}, \textbf{1}\{-Z_1 > Z_3\}\right) + \Cov\left( \textbf{1}\{Z_1 > Z_2\}, \textbf{1}\{-Z_2 > Z_3\}\right) \\
    &= \Cov\left( \textbf{1}\{Z_1 > Z_2\}, \textbf{1}\{-Z_1 > Z_3\}\right) + \Cov\left( \textbf{1}\{Z_1 < Z_2\}, \textbf{1}\{-Z_1 > Z_3\}\right) = 0\,.
\end{align*}
The only exceptions arise if $i+j=d_2$ or if $i \leq d_2$ and $j > d_2$ (or, in case of $d_2 > n,$ if $i < d_2 - n$ and $j \geq d_2 - n$), but within the moving window of pairs $(i, j) \in \mathfrak{N}_{n,d_1},$ there are only $O(d_1^2)$ such exceptional pairs. Hence, \eqref{A1} indeed dominates the asymptotic behavior of $\Cov\left(X_\inv^{B,(d_1)}, X_\des^{B,(d_2)}\right)$ in low regimes of $d_1$. \qed

\end{appendices}

\bibliographystyle{myamsplain}
\bibliography{bibliography.bib}

@book{bjorner2006combinatorics,
  title={Combinatorics of {C}oxeter groups},
  author={Björner, Anders and Brenti, Francesco},
  volume={231},
  year={2006},
  publisher={Springer Science \& Business Media}
}

@article{bona2007generalized,
  title={Generalized descents and normality},
  author={B{\'o}na, Mikl{\'o}s},
  journal={Electronic Journal of Combinatorics},
  volume={15},
  pages={N21},
  year={2008}
}

@article{bruck2019central,
  title={A central limit theorem for the two-sided descent statistic on {C}oxeter groups},
  author={Br{\"u}ck, Benjamin and R{\"o}ttger, Frank},
  journal={Electronic Journal of Combinatorics},
  volume=29,
  number=1,
  pages={P1.1},
  doi={10.37236/10744},
  year={2022}
}

@article{chang2024central,
  title={Central limit theorems for high dimensional dependent data},
  author={Chang, Jinyuan and Chen, Xiaohui and Wu, Mingcong},
  journal={Bernoulli},
  volume={30},
  number={1},
  pages={712--742},
  year={2024},
  publisher={Bernoulli Society for Mathematical Statistics and Probability}
}

@article{chatterjee2017central,
  title={A central limit theorem for a new statistic on permutations},
  author={Chatterjee, Sourav and Diaconis, Persi},
  journal={Indian Journal of Pure and Applied Mathematics},
  volume={48},
  number={4},
  pages={561--573},
  year={2017},
  publisher={Springer}
}

@article{chernozhukov2013gaussian,
  title={Gaussian approximations and multiplier bootstrap for maxima of sums of high-dimensional random vectors},
  author={Chernozhukov, Victor and Chetverikov, Denis and Kato, Kengo},
  journal={The Annals of Statistics},
  volume={41},
  number={6},
  pages={2786--2819},
  year={2013},
  publisher={Institute of Mathematical Statistics}
}

@article{chernozhukov2017central,
  title={Central limit theorems and bootstrap in high dimensions},
  author={Chernozhukov, Victor and Chetverikov, Denis and Kato, Kengo},
  journal={The Annals of Probability},
  volume={45},
  number={4},
  pages={2309--2352},
  year={2017},
  publisher={Institute of Mathematical Statistics}
}

@article {chernozhukov2023nearly,
    AUTHOR = {Chernozhukov, Victor and Chetverikov, Denis and Koike, Yuta},
     TITLE = {Nearly optimal central limit theorem and bootstrap
              approximations in high dimensions},
   JOURNAL = {The Annals of Applied Probability},
  FJOURNAL = {The Annals of Applied Probability},
    VOLUME = {33},
      YEAR = {2023},
    NUMBER = {3},
     PAGES = {2374--2425},
      ISSN = {1050-5164,2168-8737},
   MRCLASS = {60F05 (62E17)},
MRREVIEWER = {N.\ C.\ Weber},
       DOI = {10.1214/22-aap1870},
       URL = {https://doi.org/10.1214/22-aap1870},
}

@article{coxeter1935complete,
  title={The Complete Enumeration of Finite Groups of the Form {$r_i^2 = (r_ir_j)^{k_{ij}}= 1$}},
  author={Coxeter, Harold S. M.},
  journal={Journal of the London Mathematical Society},
  volume={1},
  number={1},
  pages={21--25},
  year={1935},
  publisher={Wiley Online Library}
}

@article{das2021central,
  title={Central Limit Theorem in high dimensions: The optimal bound on dimension growth rate},
  author={Das, Debraj and Lahiri, Soumendra},
  journal={Transactions of the American Mathematical Society},
  volume={374},
  number={10},
  pages={6991--7009},
  year={2021}
}

@article{de1988generalized,
  title={Generalized Eulerian numbers and the topology of the {H}essenberg variety of a matrix},
  author={de Mari, Filippo and Shayman, Mark A.},
  journal={Acta Applicandae Mathematica},
  volume={12},
  pages={213--235},
  year={1988},
  publisher={Springer}
}

@article{dorr2022extreme,
  title={Extreme Values of Permutation Statistics},
  author={D{\"o}rr, Philip and Kahle, Thomas},
  journal={The Electronic Journal of Combinatorics},
  pages={P3--10},
  year={2024}
}

@article{dorr2025joint,
  title={Joint extremes of inversions and descents of random permutations},
  author={D{\"o}rr, Philip and Heiny, Johannes},
  journal={Journal of Theoretical Probability},
  volume={38},
  number={2},
  pages={42},
  year={2025},
  publisher={Springer}
}

@article{fang2015rates,
  title={Rates of convergence for multivariate normal approximation with applications to dense graphs and doubly indexed permutation statistics},
  author={Fang, Xiao and R{\"o}llin, Adrian},
  journal={Bernoulli},
  volume={21},
  number={4},
  pages={2157--2189},
  year={2015},
  publisher={Bernoulli Society for Mathematical Statistics and Probability}
}

@article{fang2021high,
  title={{High-dimensional central limit theorems by Stein’s method}},
  author={Fang, Xiao and Koike, Yuta},
  journal={The Annals of Applied Probability},
  volume={31},
  number={4},
  pages={1660--1686},
  year={2021},
  publisher={Institute of Mathematical Statistics}
}

@article{feray2012asymptotic,
  title={Asymptotic behavior of some statistics in Ewens random permutations},
  author={F{\'e}ray, Valentin},
  journal={Discrete Mathematics \& Theoretical Computer Science},
  number={Proceedings},
  year={2012},
  publisher={Episciences. org}
}

@article {he2022central,
    AUTHOR = {He, Jimmy},
     TITLE = {A central limit theorem for descents of a {M}allows
              permutation and its inverse},
   JOURNAL = {Annales de l'Institut Henri Poincar\'{e} Probabilit\'{e}s et
              Statistiques},
  FJOURNAL = {Annales de l'Institut Henri Poincar\'{e} Probabilit\'{e}s et
              Statistiques},
    VOLUME = {58},
      YEAR = {2022},
    NUMBER = {2},
     PAGES = {667--694},
      ISSN = {0246-0203,1778-7017},
   MRCLASS = {60F05 (05E16 20F55 60B15)},
MRREVIEWER = {Chinmoy\ Bhattacharjee},
       DOI = {10.1214/21-aihp1167},
       URL = {https://doi.org/10.1214/21-aihp1167},
}

@article{kahle2020counting,
  title={Counting inversions and descents of random elements in finite {C}oxeter groups},
  author={Kahle, Thomas and Stump, Christian},
  journal={Mathematics of Computation},
  volume=89,
  number=321,
  pages={437--464},
  year=2020
}

@article{koike2021notes,
  title={Notes on the dimension dependence in high-dimensional central limit theorems for hyperrectangles},
  author={Koike, Yuta},
  journal={Japanese Journal of Statistics and Data Science},
  volume={4},
  number={1},
  pages={257--297},
  year={2021},
  publisher={Springer}
}

@article{kuchibhotla2020high,
  title={{High-dimensional CLT for Sums of Non-degenerate Random Vectors: $n^{-1/2}$-rate}},
  author={Kuchibhotla, Arun K. and Rinaldo, Alessandro},
  journal={arXiv preprint arXiv:2009.13673},
  year={2020}
}

@book{leadbetter2012extremes,
  title={Extremes and related properties of random sequences and processes},
  author={Leadbetter, Malcolm R. and Lindgren, Georg and Rootz{\'e}n, Holger},
  year={1983},
  publisher={Springer Science \& Business Media}
}

@article{meier2022central,
  title={Central limit theorems for generalized descents and generalized inversions in finite root systems},
  author={Meier, Kathrin and Stump, Christian},
  journal={Electronic Journal of Probability},
  volume={28},
  pages={1--25},
  year={2023},
  publisher={The Institute of Mathematical Statistics and the Bernoulli Society}
}

@article{pike2011convergence,
  title={Convergence rates for generalized descents},
  author={Pike, John},
  journal={Electronic Journal of Combinatorics},
  pages={P236--P236},
  year={2011}
}

@article{rottger2018asymptotics,
  title={Asymptotics of a locally dependent statistic on finite reflection groups},
  author={R{\"o}ttger, Frank},
  journal={Electronic Journal of Combinatorics},
  volume=27,
  number=2,
  pages={P2.24},
  year=2020
}

@article{sibuya1960bivariate,
  title={Bivariate extreme statistics},
  author={Sibuya, Masaaki},
  journal={Annals of the Institute of Statistical Mathematics},
  volume={11},
  number={2},
  pages={195--210},
  year={1960},
  publisher={Tokyo}
}

@book{van2000asymptotic,
  title={Asymptotic statistics},
  author={van der Vaart, Aad W.},
  volume={3},
  year={2000},
  publisher={Cambridge university press}
}

@article{sun2024central,
  title={{A Central Limit Theorem on Two-Sided Descents of Mallows Distributed Elements of Finite Coxeter Groups}},
  author={Maxwell Sun},
  journal={arXiv preprint arXiv:2406.07597},
  year={2024}
}

\end{document}